\documentclass[12pt,reqno]{amsart}
\usepackage{amsmath,amssymb,amsfonts,amscd,latexsym,amsthm,mathrsfs}
\usepackage[usenames]{color}
\usepackage[unicode]{hyperref}
\usepackage{hypbmsec}
\makeatletter
\renewcommand\thesubsection{\@arabic\c@subsection}
\makeatother

\let\wt\widetilde
\let\ol\overline
\let\le\leqslant

\let\lf\lfloor
\let\rf\rfloor
\renewcommand{\d}{{\mathrm d}}

\newcommand{\cE}{{\mathcal E}}
\newcommand{\fG}{{\mathfrak G}}
\newcommand{\fg}{{\mathfrak g}}

\newcommand{\ord}{\operatorname{ord}}
\newcommand{\perm}[1]{^{\langle#1\rangle}}

\newcommand{\R}{\mathbb{R}}

\renewcommand{\Re}{\operatorname{Re}}

\newcommand{\Res}{\operatornamewithlimits{Res}}
\newcommand{\sym}{{\operatorname{sym}}}

\newcommand{\ba}{\boldsymbol a}
\newcommand{\bb}{\boldsymbol b}
\newcommand{\be}{\boldsymbol e}
\newcommand{\bh}{\boldsymbol h}
\newcommand{\Beta}{\boldsymbol\eta}
\newcommand{\balpha}{\boldsymbol\alpha}
\newcommand{\bbeta}{\boldsymbol\beta}
\theoremstyle{remark}

\begin{document}

\title{Hypergeometric rational approximations to $\zeta(4)$}

\date{May 2019}

\author{Raffaele Marcovecchio}
\address{Dipartimento di Ingegneria e Geologia, Universit\`a di Chieti-Pescara, Viale Pindaro, 42, 65127 Pescara, Italy}
\email{raffaele.marcovecchio@unich.it}

\author{Wadim Zudilin}
\address{Department of Mathematics, IMAPP, Radboud University, PO Box 9010, 6500~GL Nijmegen, Netherlands}
\email{w.zudilin@math.ru.nl}


\begin{abstract}
We give a new hypergeometric construction of rational approximations to $\zeta(4)$,
which absorbs the earlier one from 2003 based on Bailey's ${}_9F_8$ hypergeometric integrals.
With the novel ingredients we are able to get a better control of arithmetic and produce a record irrationality measure for~$\zeta(4)$.
\end{abstract}

\maketitle

\section{Introduction}
\label{intro}

Ap\'ery's proof \cite{Ap79,Fi04,vdP79} of the irrationality of $\zeta(3)$ in the 1970s sparked research in arithmetic of the values of Riemann's zeta function $\zeta(s)$ at integers $s\ge2$.
Some particular representatives of this development include \cite{Be79,Ha90,Gu83,Ne96}, and the story culminated in a remarkable arithmetic method \cite{RV96,RV01} of Rhin and Viola to produce sharp irrationality measures for $\zeta(2)$ and $\zeta(3)$
using groups of transformations of rational approximations to the quantities.
In spite of hopes to (promptly) extend Ap\'ery's success to $\zeta(5)$ and other zeta values,
the next achievement in this direction \cite{BR01,Ri00} materialised only in the 2000s in the work of Ball and Rivoal.
The latter result helped to unify differently looking approaches for arithmetic investigations of zeta values $\zeta(s)$ and related constants under a `hypergeometric' umbrella, with some particular highlights given in \cite{Zu03,Zu04a} by one of these authors.
The hypergeometric machinery has proven to be useful in further arithmetic applications; see, for example, \cite{FSZ18,KZ19,Ma09,Zu14} for more recent achievements.

The quantity $\zeta(4)$, though known to be irrational and even transcendental, remains a natural target for testing the hypergeometry.
Ap\'ery-type approximations to the number were discovered and rediscovered on several occasions \cite{Co81,So02,Zu03} but they are not good enough to conclude about its irrationality.
In \cite{Zu03}, a general construction of rational approximations to $\zeta(4)$ is proposed, which makes use of very-well-poised hypergeometric integrals and a group of their transformations;
it leads to an estimate for the irrationality exponent of the number in question provided that a certain `denominator conjecture' for the rational approximations is valid.
The conjecture appears to be difficult enough, with its only special case established in \cite{KR07} but insufficient for arithmetic applications.
This case is usually dubbed as `most symmetric', because the group of transformations acts trivially on the corresponding approximations.

The principal goal of this work is to recast the rational approximations to $\zeta(4)$ from \cite{Zu03} in a different form (still hypergeometric!) and obtain, by these means, a better control of the arithmetic of their coefficients.
On this way we are able to produce the estimate
$$
\mu(\zeta(4))\le12.51085940\dots
$$
for the irrationality exponent of the zeta value, which is better than the conjectural one announced in \cite{Zu03}.
This is not surprising, as we do not attempt at proving the denominator conjecture from \cite{Zu03} but instead investigate the arithmetic of approximations from the different hypergeometric family.

The plan of our exposition below is as follows. In Section~\ref{sec:int} we give a Barnes-type double integral for rational approximations to $\zeta(4)$ and then, in Section~\ref{sec:sym}, work out the particular `most symmetric' case of this integral, which clearly illustrates arithmetic features of the new representation of the approximations.
We recall general settings from \cite{Zu03} in Section~\ref{sec:old} and embed the approximations into a 12-parametric family of hypergeometric-type sums that are further discussed in greater details in Section~\ref{sec:gen}.
Furthermore, Section~\ref{sec:group} reviews (and recovers) the permutation group related to the linear forms in $1$ and $\zeta(4)$ from a special subfamily of the approximations constructed.
Finally, we investigate arithmetic aspects of the general rational approximations in Section~\ref{sec:arith} and produce a calculation that leads to the new bound for $\mu(\zeta(4))$ in Section~\ref{sec:ex}.

In the text below, we intentionally avoid producing claims (in the form of propositions and lemmas) to make our exposition a storytelling rather than a traditional mathematical writing.

\section{Integral representations}
\label{sec:int}

For $k\ge2$ even, fix a \emph{generic} set of complex parameters
\begin{equation*}
\bh=(h_0,h_{-1};h_1,h_2,\dots,h_k)
\end{equation*}
satisfying the conditions
\begin{equation*}
\max\{0,\Re(h_0-h_{-1})\}<\Re h_j<\frac12\,\Re h_0
\quad\text{for}\; j=1,\dots,k,
\end{equation*}
and define as in \cite{Zu04b} the very-well-poised hypergeometric integrals
\begin{align*}
F_k'(\bh)
&=F_k'(h_0,h_{-1};h_1,h_2,\dots,h_k)
\\
&=\frac1{2\pi i}\int_{-i\infty}^{i\infty}(h_0+2t)
\frac{\prod_{j=-1}^k\Gamma(h_j+t)\cdot\Gamma(h_{-1}-h_0-t)\,\Gamma(-t)}{\prod_{j=1}^k\Gamma(1+h_0-h_j+t)}\,\d t.
\end{align*}
By Bailey's integral analogue of Dougall's theorem \cite[Section~6.6]{Ba35},
\begin{equation}
F_2'(h_0,h_{-1};h_1,h_2)
=\frac{\Gamma(h_{-1})\,\Gamma(h_1)\,\Gamma(h_2)\,\Gamma(h_1+h_{-1}-h_0)\,\Gamma(h_2+h_{-1}-h_0)}{\Gamma(1+h_0-h_1-h_2)\,\Gamma(h_1+h_2+h_{-1}-h_0)}.
\nonumber
\end{equation}
Substituting this into the iteration
\begin{align}
&
F_{k+2}'(h_0,h_{-1};h_1,\dots,h_{k-1},h_k,h_{k+1},h_{k+2})
\nonumber\\ &\quad
=\frac1{\Gamma(1+h_0-h_k-h_{k+1})\,\Gamma(1+h_0-h_k-h_{k+2})\,\Gamma(1+h_0-h_{k+1}-h_{k+2})}
\nonumber\\ &\quad\quad\times
\frac1{2\pi i}\int_{-i\infty}^{i\infty}
\Gamma(h_k+s)\,\Gamma(h_{k+1}+s)\,\Gamma(h_{k+2}+s)\,
\nonumber\\ &\quad\quad\;\times
\Gamma(1+h_0-h_k-h_{k+1}-h_{k+2}-s)
\cdot F_k'(h_0,h_{-1};-s,h_1,\dots,h_{k-1})\,\d s
\nonumber
\end{align}
obtained in \cite[Section~3]{Zu04b} we deduce consequently that
\begin{align}
&
F_4'(h_0,h_{-1};h_1,h_2,h_3,h_4)
\nonumber\\ &\quad
=\frac{\Gamma(h_{-1})\,\Gamma(h_1)\,\Gamma(h_1+h_{-1}-h_0)}{\Gamma(1+h_0-h_2-h_3)\,\Gamma(1+h_0-h_2-h_4)\,\Gamma(1+h_0-h_3-h_4)}
\nonumber\\ &\quad\quad\times
\frac1{2\pi i}\int_{-i\infty}^{i\infty}
\Gamma(h_2+s)\,\Gamma(h_3+s)\,\Gamma(h_4+s)\,\Gamma(1+h_0-h_2-h_3-h_4-s)
\nonumber\\ &\quad\quad\;\times
\frac{\Gamma(h_{-1}-h_0-s)\,\Gamma(-s)}{\Gamma(1+h_0-h_1+s)\,\Gamma(h_1+h_{-1}-h_0-s)}\,\d s
\nonumber
\end{align}
and
\begin{align}
&
F_6'(h_0,h_{-1};h_1,h_2,h_3,h_4,h_5,h_6)
\nonumber\\ &\quad
=\frac1{\Gamma(1+h_0-h_4-h_5)\,\Gamma(1+h_0-h_4-h_6)\,\Gamma(1+h_0-h_5-h_6)}
\nonumber\\ &\quad\quad\times
\frac1{2\pi i}\int_{-i\infty}^{i\infty}
\Gamma(h_4+t)\,\Gamma(h_5+t)\,\Gamma(h_6+t)\,
\nonumber\\ &\quad\quad\;\times
\Gamma(1+h_0-h_4-h_5-h_6-t)
\cdot F_4'(h_0,h_{-1};-t,h_1,h_2,h_3)\,\d t
\nonumber\displaybreak[2]\\ &\quad
=\frac1{\Gamma(1+h_0-h_4-h_5)\,\Gamma(1+h_0-h_4-h_6)\,\Gamma(1+h_0-h_5-h_6)}
\nonumber\\ &\quad\quad\times
\frac{\Gamma(h_{-1})}{\Gamma(1+h_0-h_1-h_2)\,\Gamma(1+h_0-h_1-h_3)\,\Gamma(1+h_0-h_2-h_3)}
\,\nonumber\\ &\quad\quad\times
\frac1{2\pi i}\int_{-i\infty}^{i\infty}
\Gamma(h_4+t)\,\Gamma(h_5+t)\,\Gamma(h_6+t)\,\Gamma(1+h_0-h_4-h_5-h_6-t)
\nonumber\\ &\quad\quad\times
\frac1{2\pi i}\int_{-i\infty}^{i\infty}
\Gamma(h_1+s)\,\Gamma(h_2+s)\,\Gamma(h_3+s)\,\Gamma(1+h_0-h_1-h_2-h_3-s)
\nonumber\\ &\quad\quad\;\times
\frac{\Gamma(h_{-1}-h_0-t)\,\Gamma(-t)\,\Gamma(h_{-1}-h_0-s)\,\Gamma(-s)}{\Gamma(1+h_0+s+t)\,\Gamma(h_{-1}-h_0-s-t)}\,\d s\,\d t.
\nonumber
\end{align}
Furthermore, if
$$
h_{-1}-h_0\in\mathbb Z, \quad h_0-h_1-h_2-h_3\in\mathbb Z
\quad\text{and}\quad h_0-h_4-h_5-h_6\in\mathbb Z,
$$
then the latter can be given as
\begin{align}
&
F_6'(h_0,h_{-1};h_1,h_2,h_3,h_4,h_5,h_6)
\nonumber\\ &\quad
=\frac{(-1)^{(h_{-1}-h_0)+(h_0-h_1-h_2-h_3)+(h_0-h_4-h_5-h_6)}\Gamma(h_{-1})}
{\Gamma(1+h_0-h_1-h_2)\,\Gamma(1+h_0-h_1-h_3)\,\Gamma(1+h_0-h_2-h_3)}
\nonumber\\ &\quad\quad\times
\frac1{\Gamma(1+h_0-h_4-h_5)\,\Gamma(1+h_0-h_4-h_6)\,\Gamma(1+h_0-h_5-h_6)}
\nonumber\displaybreak[2]\\ &\quad\quad\times
\frac1{2\pi i}\int_{-i\infty}^{i\infty}
\frac{\Gamma(h_4+t)\,\Gamma(h_5+t)\,\Gamma(h_6+t)}{\Gamma(1+t)\,\Gamma(1+h_0-h_{-1}+t)\,\Gamma(h_4+h_5+h_6-h_0+t)}\,\biggl(\frac\pi{\sin\pi t}\biggr)^3
\nonumber\\ &\quad\quad\;\times
\frac1{2\pi i}\int_{-i\infty}^{i\infty}
\frac{\Gamma(h_1+s)\,\Gamma(h_2+s)\,\Gamma(h_3+s)}{\Gamma(1+s)\,\Gamma(1+h_0-h_{-1}+s)\,\Gamma(h_1+h_2+h_3-h_0+s)}\,\biggl(\frac\pi{\sin\pi s}\biggr)^3
\nonumber\\ &\quad\quad\;\;\times
\frac{\Gamma(1+h_0-h_{-1}+s+t)}{\Gamma(1+h_0+s+t)}
\,\frac{\sin\pi(s+t)}\pi\,\d s\,\d t.
\label{eq5}
\end{align}

\section{The most symmetric case}
\label{sec:sym}

Equation \eqref{eq5} has an interesting structure. For example, in the most symmetric case it implies
\begin{align}
&
F_6^\sym(n)
=F_6'(3n+2,3n+2;n+1,\dots,n+1)
=\frac1{2\pi i}\int_{c_1-i\infty}^{c_1+i\infty}
\biggl(\frac{(t+1)_n}{n!}\biggr)^3\,\biggl(\frac\pi{\sin\pi t}\biggr)^3
\nonumber\\ &\quad\quad\;\times
\frac1{2\pi i}\int_{c_2-i\infty}^{c_2+i\infty}
\biggl(\frac{(s+1)_n}{n!}\biggr)^3\,\biggl(\frac\pi{\sin\pi s}\biggr)^3
\frac{(3n+1)!}{(s+t+1)_{3n+2}}\,\frac{\sin\pi(s+t)}\pi\,\d s\,\d t.
\nonumber
\end{align}
Notice that the function
$$
\frac{(3n+1)!}{(s+t+1)_{3n+2}}\,\frac{\sin\pi(s+t)}\pi
$$
is entire in both its variables, while the poles of
$$
\biggl(\frac{(s+1)_n}{n!}\biggr)^3\,\biggl(\frac\pi{\sin\pi s}\biggr)^3
$$
in a right half-plane are at $s=0,1,2,\dots$; and the latter function is analytic in the strip $-(n+1)<\Re s<0$.
A similar structure is for
$$
\biggl(\frac{(t+1)_n}{n!}\biggr)^3\,\biggl(\frac\pi{\sin\pi t}\biggr)^3.
$$
This implies that one can take $c_1,c_2\in\mathbb R$ to be any in the range $-(n+1)<c_1,c_2<0$;
we choose $c_1=c_2=c-n$ with $c=-1/3$ for our discussion below.

Now write
\[
\sin \pi (s+t) = \sin \pi s \,\cos \pi t + \cos \pi s \,\sin \pi t,
\]
so that the integral is split into the integration
\begin{align}
\frac12F_6^\sym(n)
&=\frac1{2\pi i}\int_{c-n-i\infty}^{c-n+i\infty}
\biggl(\frac{(t+1)_n}{n!}\biggr)^3\,\biggl(\frac\pi{\sin\pi t}\biggr)^3\cos\pi t
\nonumber\\ &\quad\times
\frac1{2\pi i}\int_{c-n-i\infty}^{c-n+i\infty}
\biggl(\frac{(s+1)_n}{n!}\biggr)^3\,\frac{(3n+1)!}{(s+t+1)_{3n+2}}\biggl(\frac\pi{\sin\pi s}\biggr)^2
\d s\,\d t
\label{eq7}
\end{align}
(twice, because of the symmetry $s\leftrightarrow t$).

We first deal with the internal integral in \eqref{eq7}. The rational integrand is decomposed into the sum of partial fractions:
\begin{gather*}
\biggl(\frac{(s+1)_n}{n!}\biggr)^3\,\frac{(3n+1)!}{(s+t+1)_{3n+2}}
=\sum_{k=1}^{3n+2}\frac{A_k(t)}{s+t+k},
\\
\text{where}\quad
A_k(t)=(-1)^{k-1}\binom{3n+1}{k-1}\biggl(\frac{(-t-k+1)_n}{n!}\biggr)^3
\quad\text{for}\; k=1,2,\dots,3n+2.
\end{gather*}
Then
\begin{align*}
H_n(t)
&=\frac1{2\pi i}\int_{c-n-i\infty}^{c-n+i\infty}
\biggl(\frac{(s+1)_n}{n!}\biggr)^3\,\frac{(3n+1)!}{(s+t+1)_{3n+2}}\biggl(\frac\pi{\sin\pi s}\biggr)^2\d s
\displaybreak[2]\\
&=-\sum_{\nu=-n}^\infty\frac{\partial}{\partial s}\biggl(\biggl(\frac{(s+1)_n}{n!}\biggr)^3\frac{(3n+1)!}{(s+t+1)_{3n+2}}\biggr)\bigg|_{s=\nu}
\\ \intertext{(take $\nu_0$ any from the interval $-n\le\nu_0\le0$)}
&=-\sum_{\nu=\nu_0}^\infty\frac{\partial}{\partial s}\biggl(\biggl(\frac{(s+1)_n}{n!}\biggr)^3\frac{(3n+1)!}{(s+t+1)_{3n+2}}\biggr)\bigg|_{s=\nu}
\displaybreak[2]\\
&=\sum_{\nu=\nu_0}^\infty\sum_{k=1}^{3n+2}\frac{A_k(t)}{(\nu+t+k)^2}
=\sum_{k=1}^{3n+2}A_k(t)\sum_{\nu=\nu_0}^\infty\frac1{(\nu+t+k)^2}
\displaybreak[2]\\
&=\sum_{k=1}^{3n+2}A_k(t)\biggl(\sum_{l=1}^\infty-\sum_{l=1}^{k-1}\biggr)\frac1{(t+l+\nu_0)^2}
=-\sum_{k=1}^{3n+2}A_k(t)\sum_{l=1}^{k-1}\frac1{(t+l+\nu_0)^2},
\end{align*}
because
$$
\sum_{k=1}^{3n+2}A_k(t)=0
$$
by the residue sum theorem. The choices $\nu_0=0$ and $\nu_0=-n$ lead to the equality
\begin{equation}
-\sum_{k=1}^{3n+2}A_k(t)\sum_{l=1}^{k-1}\frac1{(t+l)^2}
=H_n(t)
=-\sum_{k=1}^{3n+2}A_k(t)\sum_{l=1}^{k-1}\frac1{(t+l-n)^2};
\label{eq8}
\end{equation}
since $A_k(t)$ are polynomials, the two representations imply that the only poles of $H_n(t)$ are located at the integers
\begin{align*}
&
\{-1,-2,\dots,-3n,-(3n+1)\}\cap\{-(2n+1),-2n,\dots,n-2,n-1\}
\\ &\quad
=\{-1,-2,\dots,-2n,-(2n+1)\}.
\end{align*}
Furthermore, the function
$$
\wt H_n(t)=\biggl(\frac{(t+1)_n}{n!}\biggr)^3H_n(t)
$$
has only poles at $t=-(n+1),-(n+2),\dots,-(2n+1)$ and vanishes at $t=-1,-2,\dots,-n$.
Moreover, $\wt H_n(t)$ is in fact a rational function of degree at most $-2$ (so that it has the zero residue at infinity);
indeed, it is the sum of rational functions
$$
\biggl(\frac{(t+1)_n}{n!}\biggr)^3\frac{\partial}{\partial s}\biggl(\biggl(\frac{(s+1)_n}{n!}\biggr)^3\frac{(3n+1)!}{(s+t+1)_{3n+2}}\biggr)\bigg|_{s=\nu},
\quad\text{where}\; \nu=0,1,2,\dots,
$$
each of degree at most~$-2$ (in $t$). This means that we have a partial-fraction decomposition
$$
\wt H_n(t)
=\sum_{j=n+1}^{2n+1}\biggl(\frac{B_j}{(t+j)^2}+\frac{C_j}{t+j}\biggr)
$$
with $\sum_{j=n+1}^{2n+1}C_j=0$. With the help of the following consequence of formula~\eqref{eq8},
$$
H_n(t)
=-\sum_{l=1}^{3n+1}\frac1{(t+l-n)^2}\sum_{k=l+1}^{3n+2}A_k(t)
=-\sum_{j=-n+1}^{2n+1}\frac1{(t+j)^2}\sum_{k=j+n+1}^{3n+2}A_k(t),
$$
we find out that
\begin{align*}
B_j
&=\wt H_n(t)(t+j)^2\big|_{t=-j}
=-\biggl(\frac{(-j+1)_n}{n!}\biggr)^3\sum_{k=j+n+1}^{3n+2}A_k(-j)
\\
&=\biggl(\frac{(-j+1)_n}{n!}\biggr)^3\sum_{k=j+n+1}^{3n+2}(-1)^k\binom{3n+1}{k-1}\biggl(\frac{(j-k+1)_n}{n!}\biggr)^3
\\ \intertext{and similarly}
C_j
&=\frac{\partial}{\partial t}\bigl(\wt H_n(t)(t+j)^2\bigr)\bigg|_{t=-j}
\\
&=\frac{\partial}{\partial t}\biggl(\frac{(-t+1)_n}{n!}\biggr)^3\bigg|_{t=-j}
\cdot\sum_{k=j+n+1}^{3n+2}(-1)^k\binom{3n+1}{k-1}\biggl(\frac{(j-k+1)_n}{n!}\biggr)^3
\\ &\quad
+\biggl(\frac{(-j+1)_n}{n!}\biggr)^3\sum_{k=j+n+1}^{3n+2}(-1)^k\binom{3n+1}{k-1}
\cdot\frac{\partial}{\partial t}\biggl(\frac{(-t-k+1)_n}{n!}\biggr)^3\bigg|_{t=-j}.
\end{align*}
Note that
$$
\frac{(-j+1)_n}{n!}\in\mathbb Z, \quad \frac{(j-k+1)_n}{n!}\in\mathbb Z
$$
and
$$
d_n\cdot\frac{\partial}{\partial t}\biggl(\frac{(-t+1)_n}{n!}\biggr)\bigg|_{t=-j}\in\mathbb Z, \quad
d_n\cdot\frac{\partial}{\partial t}\biggl(\frac{(-t-k+1)_n}{n!}\biggr)\bigg|_{t=-j}\in\mathbb Z
$$
for all $j,k\in\mathbb Z$ by the standard arithmetic properties of integer-valued polynomials \cite[Lemma 4]{Zu14},
where $d_n$ denotes the least common multiple of $1,2,\dots,n$.
Furthermore, each term of the sums for $B_j$ and $C_j$ has a factor of the form
$$
\frac{(-j+1)_n}{n!}\binom{3n+1}{k-1}\frac{(j-k+1)_n}{n!}
=\binom{3n+1}{k-1}\binom{j-1}n\binom{k-j-1}n,
$$
and these quantities are all divisible by the greatest common divisor $\Phi_n$ of numbers
$$
\binom{3n+1}{a+b+1}\binom an\binom bn, \quad\text{where}\; a,b\in\mathbb Z
$$
(there are only finitely many nonzero products on the list).
Thus,
\begin{equation*}
\Phi_n^{-1}B_j\in\mathbb Z \quad\text{and}\quad \Phi_n^{-1}d_nC_j\in\mathbb Z
\quad\text{for}\; j=n+1,\dots,2n+1.
\end{equation*}

Now
\begin{align*}
\frac12F_6^\sym(n)
&=\frac1{2\pi i}\int_{c-n-i\infty}^{c-n+i\infty}\wt H_n(t)\,\biggl(\frac\pi{\sin\pi t}\biggr)^3\cos\pi t\,\d t.
\end{align*}
Since
$$
\biggl(\frac\pi{\sin\pi t}\biggr)^3\cos\pi t=\frac1{(t-\nu)^3}+O(t-\nu)
\quad\text{as}\; t\to\nu\in\mathbb Z,
$$
we have
\begin{align*}
\frac12F_6^\sym(n)
&=\sum_{\nu=-n}^\infty\Res_{t=\nu}\wt H_n(t)\,\biggl(\frac\pi{\sin\pi t}\biggr)^3\cos\pi t
=\frac12\sum_{\nu=-n}^\infty\frac{\partial^2\wt H_n(t)}{\partial t^2}\bigg|_{t=\nu}
\\
&=\frac12\sum_{\nu=-n}^\infty\sum_{j=n+1}^{2n+1}\biggl(\frac{6B_j}{(\nu+j)^4}+\frac{2C_j}{(\nu+j)^3}\biggr)
\\
&=3\sum_{j=n+1}^{2n+1}B_j\sum_{\nu=-n}^\infty\frac1{(\nu+j)^4}+\sum_{j=n+1}^{2n+1}C_j\sum_{\nu=-n}^\infty\frac1{(\nu+j)^3}
\\
&=3\sum_{j=n+1}^{2n+1}B_j\cdot\zeta(4)
-\biggl(3\sum_{j=n+1}^{2n+1}B_j\sum_{l=1}^{j-n-1}\frac1{l^4}+\sum_{j=n+1}^{2n+1}C_j\sum_{l=1}^{j-n-1}\frac1{l^3}\biggr).
\end{align*}
This implies that
$$
\frac12\,\Phi_n^{-1}d_n^4F_6^\sym(n)\in\mathbb Z\zeta(4)+\mathbb Z.
$$
In Section~\ref{sec:arith} we reveal details of computing $\Phi_n$ (and its asymptotics as $n\to\infty$); we show that $\Phi_n$ is divisible by the product over primes
\begin{equation}
\prod_{\substack{p>\sqrt{3n}\\\frac23\le\{n/p\}<1}}p.
\label{eq9}
\end{equation}
This corresponds to the `denominator conjecture' from \cite{Zu03}; for the most symmetric case in this section it was established earlier in \cite{KR07} using different hypergeometric techniques.

\section({Old approximations to \003\266(4)})%
{Old approximations to $\zeta(4)$}
\label{sec:old}

We now concentrate on a specific setting of Section~\ref{sec:int}: $k=6$ and the parameters
\begin{equation*}
\bh=(h_0,h_{-1};h_1,h_2,h_3,h_4,h_5,h_6)
\end{equation*}
are positive integers satisfying the conditions
\begin{equation*}
h_0-h_{-1}<h_j<\frac12h_0
\quad\text{for}\; j=1,2,3,4,5,6.
\end{equation*}
Define the rational function
\begin{align*}
R(t)=R(\bh;t)
&=\gamma(\bh)\,(h_0+2t)\,\frac{\prod_{j=-1}^6\Gamma(h_j+t)}{\prod_{j=-1}^6\Gamma(1+h_0-h_j+t)}
\\
&=(h_0+2t)\,
\frac{(t+1)_{h_1-1}}{(h_1-1)!}\,
\frac{(t+1+h_0-h_2)_{h_2-1}}{(h_2-1)!}\,
\\ &\qquad\times
\frac{(t+1+h_0-h_5)_{h_5+h_{-1}-h_0-1}}{(h_5+h_{-1}-h_0-1)!}\,
\frac{(t+1+h_0-h_{-1})_{h_6+h_{-1}-h_0-1}}{(h_6+h_{-1}-h_0-1)!}
\\ &\qquad\times
\frac{(h_0-h_2-h_4)!}{(t+h_2)_{h_0-h_2-h_4+1}}\,
\frac{(h_0-h_1-h_3)!}{(t+h_3)_{h_0-h_1-h_3+1}}\,
\\ &\qquad\times
\frac{(h_0-h_4-h_6)!}{(t+h_4)_{h_0-h_4-h_6+1}}\,
\frac{(h_0-h_3-h_5)!}{(t+h_5)_{h_0-h_3-h_5+1}}
\end{align*}
with
\begin{equation*}
\gamma(\bh)
=\frac{(h_0-h_2-h_4)!\,(h_0-h_1-h_3)!\,(h_0-h_4-h_6)!\,(h_0-h_3-h_5)!}
{(h_1-1)!\,(h_2-1)!\,(h_5+h_{-1}-h_0-1)!\,(h_6+h_{-1}-h_0-1)!}\,.
\end{equation*}
Then
\begin{gather*}
F(\bh)=\gamma(\bh)F_6'(\bh)=-\sum_{t=t_0}^\infty\frac{\d}{\d t}R(\bh;t)\in\mathbb Q+\mathbb Q\zeta(4)
\\
\text{with any $t_0\in\mathbb Z$},\quad
1-\min_{1\le j\le6}\{h_j\}\le t_0\le1-\max\{0,h_0-h_{-1}\},
\end{gather*}
is essentially the very-well-poised hypergeometric integral given in \cite{Zu03}; notice, however, that the arithmetic normalisation factor $\gamma(\bh)$ slightly differs from the one used in~\cite{Zu03}.
Rearranging the order of parameters in \eqref{eq5} we obtain
\begin{align}
&
F(\bh)
=\frac{(-1)^{h_{-1}+h_0+h_1+\dots+h_6}(h_{-1}-1)!\,\gamma(\bh)}
{(h_0-h_1-h_3)!\,(h_0-h_1-h_5)!\,(h_0-h_3-h_5)!}
\nonumber\\ &\quad\quad\times
\frac1{(h_0-h_2-h_4)!\,(h_0-h_2-h_6)!\,(h_0-h_4-h_6)!}
\nonumber\\ &\quad\quad\times
\frac1{2\pi i}\int_{-i\infty}^{i\infty}
\frac{\Gamma(h_2+t)\,\Gamma(h_4+t)\,\Gamma(h_6+t)}{\Gamma(1+t)\,\Gamma(1+h_0-h_{-1}+t)\,\Gamma(h_2+h_4+h_6-h_0+t)}\,\biggl(\frac\pi{\sin\pi t}\biggr)^3
\nonumber\\ &\quad\quad\;\times
\frac1{2\pi i}\int_{-i\infty}^{i\infty}
\frac{\Gamma(h_1+s)\,\Gamma(h_3+s)\,\Gamma(h_5+s)}{\Gamma(1+s)\,\Gamma(1+h_0-h_{-1}+s)\,\Gamma(h_1+h_3+h_5-h_0+s)}\,\biggl(\frac\pi{\sin\pi s}\biggr)^3
\nonumber\\ &\quad\quad\;\;\times
\frac{\Gamma(1+h_0-h_{-1}+s+t)}{\Gamma(1+h_0+s+t)}
\,\frac{\sin\pi(s+t)}\pi\,\d s\,\d t
\nonumber\displaybreak[2]\\ &\quad
=(-1)^{h_{-1}+\dots+h_6}
\nonumber\\ &\quad\quad\times
\frac1{2\pi i}\int_{-i\infty}^{i\infty}
\frac{(t+1)_{h_2-1}}{(h_2-1)!}\,
\frac{(t+h_2+h_4+h_6-h_0)_{h_0-h_2-h_6}}{(h_0-h_2-h_6)!}\,
\nonumber\\ &\quad\quad\;\;\times
\frac{(t+1+h_0-h_{-1})_{h_6+h_{-1}-h_0-1}}{(h_6+h_{-1}-h_0-1)!}\,
\biggl(\frac\pi{\sin\pi t}\biggr)^3
\nonumber\\ &\quad\quad\times
\frac1{2\pi i}\int_{-i\infty}^{i\infty}
\frac{(s+1)_{h_1-1}}{(h_1-1)!}\,
\frac{(s+h_1+h_3+h_5-h_0)_{h_0-h_1-h_5}}{(h_0-h_1-h_5)!}\,
\nonumber\\ &\quad\quad\;\;\times
\frac{(s+1+h_0-h_{-1})_{h_5+h_{-1}-h_0-1}}{(h_5+h_{-1}-h_0-1)!}\,
\biggl(\frac\pi{\sin\pi s}\biggr)^3
\nonumber\\ &\quad\quad\;\;\times
\frac{(h_{-1}-1)!}{(t+s+1+h_0-h_{-1})_{h_{-1}}}
\,\frac{\sin\pi(s+t)}\pi\,\d s\,\d t.
\label{eq10}
\end{align}
The double integral we arrive at belongs to a more general (12-parametric) family, which we are going to discuss in the next section.

\section({General approximations to \003\266(4)})%
{General approximations to $\zeta(4)$}
\label{sec:gen}

The integral in \eqref{eq10} is a special case of
\begin{align}
G(\ba,\bb)
&=\frac1{2\pi i}\int_{-i\infty}^{i\infty}
\frac{(t+b_2)_{a_2-b_2}}{(a_2-b_2)!}\,
\frac{(t+b_4)_{a_4-b_4}}{(a_4-b_4)!}\,
\frac{(t+b_6)_{a_6-b_6}}{(a_6-b_6)!}\,
\biggl(\frac\pi{\sin\pi t}\biggr)^3
\nonumber\\ &\quad\times
\frac1{2\pi i}\int_{-i\infty}^{i\infty}
\frac{(s+b_1)_{a_1-b_1}}{(a_1-b_1)!}\,
\frac{(s+b_3)_{a_3-b_3}}{(a_3-b_3)!}\,
\frac{(s+b_5)_{a_5-b_5}}{(a_5-b_5)!}\,
\biggl(\frac\pi{\sin\pi s}\biggr)^3
\nonumber\\ &\quad\;\;\times
\frac{(b_0-a_0-1)!}{(t+s+a_0)_{b_0-a_0}}
\,\frac{\sin\pi(s+t)}\pi\,\d s\,\d t,
\label{eq11}
\end{align}
where the integral parameters
\begin{equation}
\ba=(a_0;a_1,a_2,a_3,a_4,a_5,a_6) \quad\text{and}\quad
\bb=(b_0;b_1,b_2,b_3,b_4,b_5,b_6)
\label{eq12}
\end{equation}
are subject to the conditions
\begin{equation}
\begin{aligned}
b_0-a_0-2
&\ge(a_1+a_3+a_5)-(b_1+b_3+b_5),
\\
b_0-a_0-2
&\ge(a_2+a_4+a_6)-(b_2+b_4+b_6),
\end{aligned}
\label{eq13}
\end{equation}
and
\begin{equation}
\max\{b_1,b_3,b_5\}\le\min\{a_1,a_3,a_5\}, \quad
\max\{b_2,b_4,b_6\}\le\min\{a_2,a_4,a_6\}.
\nonumber
\end{equation}
Note that simultaneous shifts
of $a_0,a_1,a_3,a_5$ and $b_0,b_1,b_3,b_5$ by the same integer does not affect $G(\ba,\bb)$;
the same is true for simultaneous shifts of $a_0,a_2,a_4,a_6$ and $b_0,b_2,b_4,b_6$.
(In particular, the shifts by given $1-b_1$ and $1-b_2$, respectively, allow to assume that $b_1=b_2=1$.)
The latter two symmetries potentially leave 12 out of 14 parameters \eqref{eq12} independent.
Furthermore, we choose
\begin{equation}
\ba^*=(a_0;a_1^*,a_2^*,a_3^*,a_4^*,a_5^*,a_6^*) \quad\text{and}\quad
\bb^*=(b_0;b_1^*,b_2^*,b_3^*,b_4^*,b_5^*,b_6^*)
\label{eq15}
\end{equation}
to be a reordering of the parameters \eqref{eq12} (so that \eqref{eq15} and \eqref{eq12} coincide as multi-sets) such that
\begin{equation}
a_1^*\le a_3^*\le a_5^*, \quad b_1^*\le b_3^*\le b_5^*
\quad\text{and}\quad
a_2^*\le a_4^*\le a_6^*, \quad b_2^*\le b_4^*\le b_6^*.
\nonumber
\end{equation}
Additionally, we assume
\begin{equation}
a_0+1\ge b_3^*+b_4^*.
\label{eq17}
\end{equation}

Similarly to the most symmetric case in Section~\ref{sec:sym}, we may choose the integration paths in \eqref{eq11} 
to be the vertical lines $\{c_1+iy:y\in\R\}$ for $s$ and $\{c_2+iy:y\in\R\}$ for~$t$, 
with
\[
-a_1^*<c_1<1-b_5^*, \quad -a_2^*<c_2<1-b_6^*,
\]
and we take $c_1=1/3-a_1^*$ and $c_2=1/3-a_2^*$. Also, the 
rational function in $s$ and $t$ at the integrand in \eqref{eq11} has degree at most  $-2$ 
both in $s$ and in $t$, and the functions
\[
\frac{1}{\sin \pi s}, \quad \frac{\cos \pi s}{(\sin \pi s)^2}\quad 
		\text{and} \quad \frac{1}{\sin \pi t}
\]
are bounded in their respective integration domains. By
\[
\sin \pi (s+t) = \sin \pi s \,\cos \pi t + \cos \pi s\, \sin \pi t
\]
the integral $G(\ba,\bb)$ is  split into two absolutely convergent 
integrals, and, after interchanging the order of integrations in $s$ and 
in $t$ in the second integral, we obtain 
\begin{align}
G(\ba,\bb)
&=\frac1{2\pi i}\int_{-i\infty}^{i\infty}
\frac{(t+b_2)_{a_2-b_2}}{(a_2-b_2)!}\,
\frac{(t+b_4)_{a_4-b_4}}{(a_4-b_4)!}\,
\frac{(t+b_6)_{a_6-b_6}}{(a_6-b_6)!}\,
\biggl(\frac\pi{\sin\pi t}\biggr)^3 \cos\pi t
\nonumber\\ &\quad\times
\frac1{2\pi i}\int_{-i\infty}^{i\infty}
\frac{(s+b_1)_{a_1-b_1}}{(a_1-b_1)!}\,
\frac{(s+b_3)_{a_3-b_3}}{(a_3-b_3)!}\,
\frac{(s+b_5)_{a_5-b_5}}{(a_5-b_5)!}\,
\biggl(\frac\pi{\sin\pi s}\biggr)^2
\nonumber\\ &\quad\;\;\times
\frac{(b_0-a_0-1)!}{(t+s+a_0)_{b_0-a_0}} \,\d s\,\d t 		
\nonumber\\ &\; + \text{a similar integral with } a_j,b_j 
					\text{ changed to } a_{7-j},b_{7-j} \text{ for } j=1,\dots,6.
\label{eq18}
\end{align}
As already seen in the most symmetric case, the integral
\begin{align}
H(t)
&=H(a_0,a_1,a_3,a_5;b_0,b_1,b_3,b_5;t)
\nonumber\\
&=\frac{1}{2\pi i}
\int_{c_1-i\infty}^{c_1+i\infty} \frac{(s+b_1)_{a_1-b_1}}{(a_1-b_1)!}\,
\frac{(s+b_3)_{a_3-b_3}}{(a_3-b_3)!}\,
\frac{(s+b_5)_{a_5-b_5}}{(a_5-b_5)!}\,
\nonumber\\ &\quad\times
\frac{(b_0-a_0-1)!}{(t+s+a_0)_{b_0-a_0}}
\biggl(\frac\pi{\sin\pi s}\biggr)^2 \,\d s
\label{eq19}
\end{align}
is a rational function in $t$, and we may even vary $c_1$ in 
the interval $-a_3^*<c_1<1-b_3^*$, because a power of $\sin \pi s$ is dropped 
in the denominator of \eqref{eq19} with respect to the integral \eqref{eq11}.  
In executing this, we do not have to take care of possible poles coming from 
$(t+s+a_0)_{b_0-a_0}$, because it never vanishes if $t$ is chosen in an appropriate 
region of the complex plane, and two rational functions that coincide in such 
a region must coincide everywhere.

Explicitly, we have
\[
\frac{(s+b_1)_{a_1-b_1}}{(a_1-b_1)!}\,
\frac{(s+b_3)_{a_3-b_3}}{(a_3-b_3)!}\,
\frac{(s+b_5)_{a_5-b_5}}{(a_5-b_5)!}\,
\frac{(b_0-a_0-1)!}{(t+s+a_0)_{b_0-a_0}}
= \sum_{k=a_0}^{b_0-1} \frac{A_k(t)}{t+s+k}, 
\]
where
\begin{align}
A_k(t)
&= (-1)^{k+a_0} \binom{b_0-a_0-1}{k-a_0} 
\frac{(-t-k+b_1)_{a_1-b_1}}{(a_1-b_1)!}
\nonumber\\ &\quad\times
\frac{(-t-k+b_3)_{a_3-b_3}}{(a_3-b_3)!}
\frac{(-t-k+b_3)_{a_5-b_5}}{(a_5-b_5)!}
\quad\text{for}\; k=a_0,\dots,b_0-1
\label{eq20}
\end{align}
satisfy $\sum_{k=a_0}^{b_0-1} A_k(t)=0$. Then
\begin{align}
H(t)
=&-\sum_{\nu=\nu_0}^\infty\frac{\partial}{\partial s}\sum_{k=a_0}^{b_0-1} 
\frac{A_k(t)}{t+s+k}\bigg|_{s=\nu}
=\sum_{k=a_0}^{b_0-1}A_k(t) \sum_{\nu=\nu_0}^\infty\frac1{(\nu+t+k)^2}
\nonumber\\
=&-\sum_{k=a_0}^{b_0-1}A_k(t) \sum_{l=a_0}^{k-1}\frac1{(t+l+\nu_0)^2},
\label{eq21}
\end{align}
where $\nu_0$ is any integer in the interval $1-a_3^*\le\nu_0\le 1-b_3^*$. 
Since all $A_k(t)$ are polynomials, the poles of function \eqref{eq21} 
are only possible at
$$
t=a_3^*-b_0+1,a_3^*-b_0+2,\dots,b_3^*-a_0-1.
$$
For a similar reason, with $\nu_0$ in the larger interval $1-a_5^*\le\nu_0\le 1-b_1^*$, the function
\[
I(t)=\sum_{\nu=\nu_0}^\infty\sum_{k=a_0}^{b_0-1} \frac{A_k(t)}{t+s+k}\bigg|_{s=\nu}
=-\sum_{k=a_0}^{b_0-1}A_k(t) \sum_{l=a_0}^{k-1}\frac1{t+l+\nu_0}
\]
has only poles possible at 
\[
t=a_5^*-b_0+1,a_5^*-b_0+2,\dots,b_1^*-a_0-1.
\]
Since
\[
H(t) (t+l+\nu_0)^2\big|_{t=-l-\nu_0}
=\sum_{k=l+1}^{b_0} A_k(-l-\nu_0)
=I(t) (t+l+\nu_0)\big|_{t=-l-\nu_0} 
\]
when $1-a_3^*\le\nu_0\le 1-b_3^*$, it follows that the set of double poles of $H(t)$ 
coincides with the set of simple poles of $I(t)$, and therefore is also contained at integers
in $[a_5^*-b_0+1,b_1^*-a_0-1]$; however, $H(t)$ may still possess simple poles at integers
in $[a_3^*-b_0+1,b_3^*-a_0-1]$.
Arguing as in Section~\ref{sec:sym} we arrive at the partial-fraction decomposition
\begin{align}
\wt H(t)
=\wt H(\ba,\bb;t)
&=\frac{(t+b_2)_{a_2-b_2}}{(a_2-b_2)!}\,
\frac{(t+b_4)_{a_4-b_4}}{(a_4-b_4)!}\,
\frac{(t+b_6)_{a_6-b_6}}{(a_6-b_6)!}\,
H(t)
\nonumber\\
&=\sum_{j=1+a_0-b_1^*}^{b_0-a_5^*-1}\frac{B_j}{(t+j)^2}
+\sum_{j=1+a_0-b_3^*}^{b_0-a_3^*-1}\frac{C_j}{t+j},
\label{eq22}
\end{align}
because the rational function $\wt H(t)$ has degree at most $-2$ by~\eqref{eq13}.
Noticing that the expression
$$
\frac{(t+b_2)_{a_2-b_2}}{(a_2-b_2)!}\,
\frac{(t+b_4)_{a_4-b_4}}{(a_4-b_4)!}\,
\frac{(t+b_6)_{a_6-b_6}}{(a_6-b_6)!}
$$
has at least simple zeroes at $t=1-a_6^*,\dots,-b_2^*$ 
and at least double zeroes at $t=1-a_4^*,2-a_4^*,\dots,-b_4^*$ and taking into 
account condition \eqref{eq17}, we find out that $\wt H(t)$ does not have poles in the half-plane $\Re t>c_2$,
hence  the expansion \eqref{eq22} `shortens' to
\begin{align*}
\wt H(t)
=\sum_{j=a_2^*}^{b_0-a_5^*-1}\frac{B_j}{(t+j)^2}
+\sum_{j=a_2^*}^{b_0-a_3^*-1}\frac{C_j}{t+j}.
\end{align*}
In fact, the second sum is over the interval $\max\{a_4^*,1+a_0-b_3^*\}\le j\le b_0-a_3^*-1$,
while the first one is over $\max\{a_6^*,1+a_0-b_1^*\}\le j\le b_0-a_5^*-1$ and may be even empty if the interval is empty.
With the explicit expressions \eqref{eq20} and \eqref{eq21} (used, for example, with $\nu_0=1-a_3^*$) 
in mind, we conclude that the coefficients
$$
B_j=\wt H(t)(t+j)^2
\quad\text{and}\quad
C_j=\frac{\partial}{\partial t}\bigl(\wt H(t)(t+j)^2\bigr)\Big|_{t=-j}
\quad\text{for}\; j\in\mathbb Z
$$
satisfy
\begin{gather}
B_j\in\mathbb Z, \quad d_mC_j\in\mathbb Z
\nonumber
\\
\text{with}\; m=\max\{a_1-b_1,a_2-b_2,a_3-b_3,a_4-b_4,a_5-b_5,a_6-b_6\},
\nonumber
\end{gather}
but also
\begin{align}
\ord_pB_j, \; \ord_p(d_mC_j)
&\ge\min_{j,k\in\mathbb Z}
\biggl(\biggl\lf\frac{b_0-a_0-1}p\biggr\rf
-\biggl\lf\frac{k-a_0}p\biggr\rf
-\biggl\lf\frac{b_0-k-1}p\biggr\rf
\nonumber\\ &\qquad
+\sum_{r\in\{2,4,6\}}\biggl(\biggl\lf\frac{j-b_r}p\biggr\rf-\biggl\lf\frac{j-a_r}p\biggr\rf-\biggl\lf\frac{a_r-b_r}p\biggr\rf\biggr)
\nonumber\\ &\qquad
+\sum_{r\in\{1,3,5\}}\biggl(\biggl\lf\frac{k-j-b_r}p\biggr\rf-\biggl\lf\frac{k-j-a_r}p\biggr\rf-\biggl\lf\frac{a_r-b_r}p\biggr\rf\biggr)
\biggr)
\nonumber
\end{align}
for primes $p>\sqrt{b_0-a_0}$ (see \cite[Lemmas 17, 18]{Zu04a}).
Furthermore,
\begin{align}
&
\frac1{2\pi i}\int_{c_2-i\infty}^{c_2+i\infty}\wt H(t)\biggl(\frac\pi{\sin\pi t}\biggr)^3\cos\pi t\, \d t
\nonumber\\ &\quad
=\sum_{\nu=1-a_2^*}^\infty\biggl(\sum_{j=a_2^*}^{b_0-a_5^*-1}\frac{3B_j}{(\nu+j)^4}
+\sum_{j=a_2^*}^{b_0-a_3^*-1}\frac{C_j}{(\nu+j)^3}\biggr)
\nonumber\\ &\quad
=3\zeta(4)\sum_{j=\max\{a_6^*,1+a_0-b_1^*\}}^{b_0-a_5^*-1}B_j
-\biggl(3\sum_{j=a_2^*}^{b_0-a_5^*-1}B_j\sum_{l=1}^{j-a_2^*}\frac1{l^4}
+\sum_{j=a_2^*}^{b_0-a_3^*-1}C_j\sum_{l=1}^{j-a_2^*}\frac1{l^3}\biggr),
\nonumber
\end{align}
where $\sum_jC_j=0$ was implemented.
Performing the same way for the second double integral in~\eqref{eq18} we conclude that
\begin{gather}
G(\ba,\bb)=B(\ba,\bb)\zeta(4)-C(\ba,\bb),
\quad\text{where}\;
B\in\mathbb Z, \; d_{m_1}^3d_{m_2}C\in\mathbb Z 
\label{eq26}
\\
\text{with}\; m_1=\max\{b_0-a_2^*-a_3^*-1,b_0-a_1^*-a_4^*-1\},
\nonumber
\\
\; m_2=\max\{b_0-a_2^*-a_5^*-1,b_0-a_1^*-a_6^*-1,a_1-b_1,\dots,a_6-b_6\},
\nonumber
\end{gather}
and
\begin{align}
\ord_pB, \; 4+\ord_pC
&\ge\min_{j,l\in\mathbb Z}
\biggl(\biggl\lf\frac{b_0-a_0-1}p\biggr\rf
-\biggl\lf\frac{j+l-a_0}p\biggr\rf
-\biggl\lf\frac{b_0-j-l-1}p\biggr\rf
\nonumber\\ &\qquad
+\sum_{r\in\{2,4,6\}}\biggl(\biggl\lf\frac{j-b_r}p\biggr\rf-\biggl\lf\frac{j-a_r}p\biggr\rf-\biggl\lf\frac{a_r-b_r}p\biggr\rf\biggr)
\nonumber\\ &\qquad
+\sum_{r\in\{1,3,5\}}\biggl(\biggl\lf\frac{l-b_r}p\biggr\rf-\biggl\lf\frac{l-a_r}p\biggr\rf-\biggl\lf\frac{a_r-b_r}p\biggr\rf\biggr)
\biggr)
\label{eq27}
\end{align}
for primes $p>\sqrt{b_0-a_0-2}$.

\medskip
Finally, we remark that condition \eqref{eq17} is conventional 
(and happens to hold in our applications, even in the form of equality $b_3^*+b_4^*=a_0-1$) but can be potentially dropped without 
significant arithmetic losses. For example, if $b_1^*+b_2^*>a_0-1$ then 
the partial-fraction decomposition \eqref{eq22} translates into
\begin{align*}
\wt H(t)
=\sum_{j=1+a_0-b_1^*}^{b_2^*-1}\frac{B_j}{(t+j)^2}+\sum_{j=a_6^*}^{b_0-a_5^*-1}\frac{B_j}{(t+j)^2}
+\sum_{j=1+a_0-b_3^*}^{b_4^*-1}\frac{C_j}{t+j}+\sum_{j=a_4^*}^{b_0-a_3^*-1}\frac{C_j}{t+j},
\end{align*}
so that there are poles of $\wt H(t)$ to the right of the contour $\Re t=c_2$. The corresponding residues of the integrand are
\begin{gather*}
\Res_{t=-j}\wt H(t)\biggl(\frac\pi{\sin\pi t}\biggr)^3\cos\pi t
=D_j-6\zeta(4)\,B_j
\\
\quad\text{with}\; 
D_j=\frac1{24}\,\frac{\partial^4}{\partial t^4}\bigl(\wt H(t)(t+j)^2\bigr)\Big|_{t=-j}
=\frac12\,\frac{\partial^2}{\partial t^2}\biggl(\wt H(t)-\frac{B_j}{(t+j)^2}-\frac{C_j}{t+j}\biggr)\bigg|_{t=-j},
\end{gather*}
where $j$ is an integer in the interval $1+a_0-b_3^*\le j\le b_4^*-1$ and we use the expansion
\begin{equation*}
\biggl(\frac\pi{\sin\pi t}\biggr)^3\cos\pi t
=\frac1{(t+j)^3}-6\zeta(4)\,(t+j)+O\bigl((t+j)^3\bigr)
\quad\text{as}\; t\to-j.
\end{equation*}
Proceeding as above we deduce that
\begin{align}
&
\frac1{2\pi i}\int_{c_2-i\infty}^{c_2+i\infty}\wt H(t)\biggl(\frac\pi{\sin\pi t}\biggr)^3\cos\pi t\, \d t
=\sum_{\nu=1-a_2^*}^\infty\Res_{t=\nu}\wt H(t)\biggl(\frac\pi{\sin\pi t}\biggr)^3\cos\pi t
\nonumber\\ &\;
=-6\zeta(4)\sum_{j=1+a_0-b_1^*}^{b_2^*-1}B_j
+3\sum_{j=1+a_0-b_1^*}^{b_2^*-1}B_j\sum_{\substack{l=j+1-a_2^*\\l\ne0}}^\infty\frac1{l^4}
+3\sum_{j=a_6^*}^{b_0-a_5^*-1}B_j\sum_{l=j+1-a_2^*}^\infty\frac1{l^4}
\nonumber\\ &\;\quad
+\sum_{j=1+a_0-b_3^*}^{b_4^*-1}C_j\sum_{\substack{l=j+1-a_2^*\\l\ne0}}^\infty\frac1{l^3}
+\sum_{j=a_4^*}^{b_0-a_3^*-1}C_j\sum_{l=j+1-a_2^*}^\infty\frac1{l^3},
\nonumber
\end{align}
which is again seen to be a linear form in $\mathbb Z\zeta(4)+\mathbb Q$.

\section({The group structure for \003\266(4)})%
{The group structure for $\zeta(4)$}
\label{sec:group}

Following \cite{Zu03}, to any set of parameters $\bh$ from Section~\ref{sec:old} 
we assign the 27-element multiset of \emph{nonnegative} integers
\begin{equation}
\begin{gathered}
e_{0j}=h_j-1, \; \ol e_{0j}=h_j+h_{-1}-h_0-1 \quad\text{for}\; 1\le j\le 6,
\\
e_{jk}=h_0-h_j-h_k \quad\text{for}\; 1\le j<k\le 6,
\end{gathered}
\label{eq28}
\end{equation}
and set $H(\be)=F(\bh)$ for the quantity defined in that section. By the construction,
$$
\gamma(\bh)^{-1}F(\bh)=\frac{e_{01}!\,e_{02}!\,\ol e_{05}!\,\ol e_{06}!}{e_{13}!\,e_{24}!\,e_{35}!\,e_{46}!}\,H(\be)
$$
is invariant under any permutation of the parameters $h_1,h_2,\dots,h_6$ 
(which we can view as the `$\bh$-trivial' action). Clearly, any such permutation 
induces the corresponding permutation of the parameter set~\eqref{eq28}.

On the other hand, it is seen from \eqref{eq11} that the quantity
$$
\prod_{j=1}^6(a_j-b_j)!\cdot G(\ba,\bb)
$$
does not change when the parameters in either collection $a_1,a_3,a_5$ 
or $a_2,a_4,a_6$ permute; we can regard such permutations as `$\ba$-trivial'. 
(The same effect is produced by `$\bb$-trivial' permutations, when we change 
the order in $b_1,b_3,b_5$ or $b_2,b_4,b_6$.) We can also add to the list 
the `trivial' involution
$$
\mathfrak{i}\colon a_j\leftrightarrow a_{7-j}, \; b_j\leftrightarrow b_{7-j}
\quad\text{for}\; j=1,\dots,6,
$$
which reflects the symmetry $s\leftrightarrow t$ 
of the double integral \eqref{eq11}. In addition, we recall that $G(\ba,\bb)$ 
is left unchanged by the simultaneous shifts of $a_0,a_1,a_3,a_5$ and $b_0,b_1,b_3,b_5$ 
(or of $a_0,a_2,a_4,a_6$ and $b_0,b_2,b_4,b_6$, respectively) by the same integer. 
We regard the action of all these transformations (permutations, shifts and involution) 
and their compositions as the `$(\ba,\bb)$-trivial' action. 

By setting
\begin{equation}
\begin{gathered}
a_0=1+h_0-h_{-1}, \quad a_j=h_j \quad\text{for}\; j=1,\dots,6,
\\
b_0=1+h_0, \quad b_1=b_2=1, \quad b_5=b_6=1+h_0-h_{-1},
\\
b_3=h_1+h_3+h_5-h_0, \quad b_4=h_2+h_4+h_6-h_0.
\end{gathered}
\label{eq29}
\end{equation}
we have $F(\bh)=G(\ba,\bb)$. If we request the condition
\begin{equation}
h_{-1}=2+3h_0-(h_1+h_2+h_3+h_4+h_5+h_6)
\label{eq30}
\end{equation}
to hold, then the shift of $h_0,h_1,h_3,h_5$ 
by $1+h_0-h_1-h_3-h_5$, that is, the transformation
\begin{multline*}
\mathfrak{b}_{135}\colon\bh \mapsto (1+2h_0-h_1-h_3-h_5,h_{-1};1+h_0-h_3-h_5,h_2, \\
			1+h_0-h_1-h_5,h_4,1+h_0-h_1-h_3,h_6), 
\end{multline*}
induces the composition of the shift of $a_0,a_1,a_3,a_5$ and $b_0,b_1,b_3,b_5$ 
by $1+h_0-h_1-h_3-h_5$ and the permutation $(b_1\ b_3)(b_4\ b_6)$. Therefore 
$\mathfrak{b}_{135}$, which also induces the permutation
\[
\mathfrak{b}=(e_{01}\ e_{35})(e_{03}\ e_{15})(e_{05}\ e_{13})
			(\ol e_{02}\ e_{46})(\ol e_{04}\ e_{26})(\ol e_{06}\ e_{24})
\]
on the parameter set \eqref{eq28},
is an $(\ba,\bb)$-trivial transformation. As a consequence, the quantity
\[
\prod_{j=1}^6(a_j-b_j)!\cdot G(\ba,\bb) = e_{01}!\,e_{02}!\,\ol e_{05}!\,\ol e_{06}!\,
											e_{15}!\,e_{26}!\,H(\be)
\]
does not change by the action of the permutation $\mathfrak{b}$. 
We remark that \eqref{eq30} is a very natural condition for the application 
of the $(\ba,\bb)$-trivial action to $F(\bh)$. Indeed, by \eqref{eq29} we have 
$b_1-b_5=b_2-b_6$, and \eqref{eq30} is equivalent to $b_1-b_3=b_4-b_6$, or to 
$b_3-b_5=b_2-b_4$. 

Taking the multiset
\begin{equation}
\cE=\{e_{03},e_{04},e_{05},e_{06},\ol e_{01},\ol e_{02},\ol e_{03},\ol e_{04},e_{13},e_{24},e_{35},e_{46}\}
\nonumber
\end{equation}
we conclude that the quantity
\begin{equation}
\frac{H(\be)}{\prod_{e\in\cE}e!}
=\frac1{\prod_{j=1}^6e_{0j}!\,\ol e_{0j}!}
\cdot\frac{e_{01}!\,e_{02}!\,\ol e_{05}!\,\ol e_{06}!}{e_{13}!\,e_{24}!\,e_{35}!\,e_{46}!}\,H(\be)
\nonumber
\end{equation}
is invariant under the $\bh$-trivial permutations and, by
\begin{equation}
\frac{H(\be)}{\prod_{e\in\cE}e!}
=\frac{\prod_{j=1}^6(a_j-b_j)!\cdot G(\ba,\bb)}{e_{13}!\,e_{15}!\,e_{35}!\,e_{24}!\,e_{26}!\,e_{46}!\,\prod_{j=1}^6e_{0j}!\,\ol e_{0j}!},
\nonumber
\end{equation}
also under the permutation $\mathfrak{b}$.

The permutation group of the multiset \eqref{eq28}, which is generated by 
all $\bh$-trivial and the permutation $\mathfrak{b}$, coincides with 
the group $\fG$ (of order $51840$) considered in \cite{Zu03}.
(Note that the group contains the above involution $\mathfrak{i}$ as well.)
By these means we also recover the invariance of the quantity
\begin{equation}
\frac{H(\be)}{\Pi(\be)},
\quad\text{where}\;\;
\Pi(\be)=e_{03}!\,e_{04}!\,e_{05}!\,e_{06}!\,\ol e_{01}!\,\ol e_{02}!\,\ol e_{03}!\,\ol e_{04}!\,e_{13}!\,e_{24}!\,e_{35}!\,e_{46}!\,,
\nonumber
\end{equation}
under the action of $\fG$ and corresponding to the arithmetic normalisation 
of $H(\be)=F(\bh)=G(\ba,\bb)$ in Section~\ref{sec:old}. 

Because our access to the arithmetic of coefficients of linear 
forms $H(\be)\in\mathbb Z\zeta(4)+\mathbb Q$ is performed through 
their $G(\ba,\bb)$-representation, we will be interested in collecting 
a set of representatives which are distinct modulo $(\ba,\bb)$-trivial 
transformations. For a generic set of integral parameters $\bh$ subject 
to \eqref{eq30}, such set of representatives contains 120 different elements. 
Indeed, by \eqref{eq29} and \eqref{eq30} the subgroup of all the $(\ba,\bb)$-trivial 
permutations in $\fG$ contains $3!^3\cdot 2!=432$ elements, and is generated by:
\begin{itemize}
\item the $\ba$- and $\bh$-trivial permutations $(h_1\ h_3)$ and $(h_3\ h_5)$;
\item the $\ba$- and $\bh$-trivial permutations $(h_2\ h_4)$ and $(h_4\ h_6)$;
\item the $\bb$-trivial permutation $(b_1\ b_3)(b_4\ b_6)$ (that is, by $\mathfrak{b}_{135}$);
and
\item the involution $\mathfrak{i}$ (that is, by $(h_1\ h_6)(h_2\ h_5)(h_3\ h_4)$).
\end{itemize}
This subgroup also contains $(b_2\ b_4)(b_3\ b_5)$ 
(namely, $\mathfrak{b}_{246}=\mathfrak{i}\mathfrak{b}_{135}\mathfrak{i}$), 
and is isomorphic to $\mathfrak{S}_3^3\times\mathfrak{S}_2$.
Now, the group $\fG$ is generated by $(h_1\ h_3)$, $(h_3\ h_5)$, $(h_2\ h_4)$, $(h_4\ h_6)$, $\mathfrak{b}_{135}$  and $(h_3\ h_4)$.
Note that $(h_1\ h_3)$, $(h_3\ h_5)$, $(h_2\ h_4)$  and $(h_4\ h_6)$ commute with $\mathfrak{b}_{135}$,
while $(h_3\ h_4)$ acts on $(\ba,\bb)$ by $a_3	\leftrightarrow a_4$, $b_3\mapsto b_3+a_4-a_3$, 
$b_4\mapsto b_4+a_3-a_4$ (and leaves $a_i, b_i$ unchanged for $i\not=3,4$).
Hence there are exactly $|\fG|/432=120$ elements in $\fG$ that are distinct 
modulo the $(\ba,\bb)$-trivial subgroup, each for any simultaneous choice of a 
subset $\{a_1,a_3,a_5\}$ (or $\{a_2,a_4,a_6\}$) of $\{h_1,\dots,h_6\}$ (among 
all $\binom{6}{3}=20$ such subsets) and of a permutation in the $\bb$-trivial subgroup 
(of $3!=6$ elements) generated by $\mathfrak{b}_{135}$ and $\mathfrak{b}_{246}$.

\section{Arithmetic of linear forms}
\label{sec:arith}

In order to compute the minimum on the right-hand side of \eqref{eq27}, we distinguish two different situations:
(a) $j+l-a_0$ is coprime with~$p$, and (b) $j+l-a_0$ is divisible by~$p$. In case~(a), we get $\lf(j+l-a_0)/p\rf=\lf(j+l-a_0-1)/p\rf$, so that the minimum in \eqref{eq27} is greater or equal than
\begin{align}
\Omega_1(\ba,\bb;p)
&=\min_{j,l\in\mathbb Z}
\biggl(\biggl\lf\frac{b_0-a_0-2}p\biggr\rf
-\biggl\lf\frac{j+l-a_0-1}p\biggr\rf
-\biggl\lf\frac{b_0-j-l-1}p\biggr\rf
\nonumber\\ &\qquad
+\sum_{r\in\{2,4,6\}}\biggl(\biggl\lf\frac{j-b_r}p\biggr\rf-\biggl\lf\frac{j-a_r}p\biggr\rf-\biggl\lf\frac{a_r-b_r}p\biggr\rf\biggr)
\nonumber\\ &\qquad
+\sum_{r\in\{1,3,5\}}\biggl(\biggl\lf\frac{l-b_r}p\biggr\rf-\biggl\lf\frac{l-a_r}p\biggr\rf-\biggl\lf\frac{a_r-b_r}p\biggr\rf\biggr)
\biggr).
\label{eq128}
\end{align}
In case (b), we have $l=-j+a_0+\mu p$ for some $\mu\in\mathbb Z$ and
\begin{align}
&
\biggl\lf\frac{b_0-a_0-1}p\biggr\rf
-\biggl\lf\frac{j+l-a_0}p\biggr\rf
-\biggl\lf\frac{b_0-j-l-1}p\biggr\rf
\nonumber\\ &\qquad
+\sum_{r\in\{2,4,6\}}\biggl(\biggl\lf\frac{j-b_r}p\biggr\rf-\biggl\lf\frac{j-a_r}p\biggr\rf-\biggl\lf\frac{a_r-b_r}p\biggr\rf\biggr)
\nonumber\\ &\qquad
+\sum_{r\in\{1,3,5\}}\biggl(\biggl\lf\frac{l-b_r}p\biggr\rf-\biggl\lf\frac{l-a_r}p\biggr\rf-\biggl\lf\frac{a_r-b_r}p\biggr\rf\biggr)
\nonumber\displaybreak[2]\\ &\quad
=\sum_{r\in\{2,4,6\}}\biggl(\biggl\lf\frac{j-b_r}p\biggr\rf-\biggl\lf\frac{j-a_r}p\biggr\rf-\biggl\lf\frac{a_r-b_r}p\biggr\rf\biggr)
\nonumber\\ &\quad\qquad
+\sum_{r\in\{1,3,5\}}\biggl(\biggl\lf\frac{-j+a_0-b_r}p\biggr\rf-\biggl\lf\frac{-j+a_0-a_r}p\biggr\rf-\biggl\lf\frac{a_r-b_r}p\biggr\rf\biggr)
\nonumber\displaybreak[2]\\ &\quad
=\sum_{r\in\{2,4,6\}}\biggl(\biggl\lf\frac{j-b_r}p\biggr\rf-\biggl\lf\frac{j-a_r}p\biggr\rf-\biggl\lf\frac{a_r-b_r}p\biggr\rf\biggr)
\nonumber\\ &\quad\qquad
+\sum_{r\in\{1,3,5\}}\biggl(\biggl\lf\frac{j+a_r-a_0-1}p\biggr\rf-\biggl\lf\frac{j+b_r-a_0-1}p\biggr\rf-\biggl\lf\frac{a_r-b_r}p\biggr\rf\biggr)
\nonumber
\end{align}
for primes $p>\sqrt{b_0-a_0-2}$, where the property $\lf\alpha+\mu\rf=\lf\alpha\rf+\mu$ was used, together with the passage
$$
\biggl\{\frac{a-1}p\biggr\}+\biggl\{\frac{-a}p\biggr\}=\frac{p-1}p,
\quad\text{where}\; a\in\mathbb Z,
$$
for the fractional part $\{\alpha\}=\alpha-\lf\alpha\rf$ of a number. This means that in case (b) the minimum in~\eqref{eq27} is equal to
\begin{align}
\Omega_2(\ba,\bb;p)
&=\min_{j\in\mathbb Z}
\biggl(\sum_{r\in\{2,4,6\}}\biggl(\biggl\lf\frac{j-b_r}p\biggr\rf-\biggl\lf\frac{j-a_r}p\biggr\rf-\biggl\lf\frac{a_r-b_r}p\biggr\rf\biggr)
\nonumber\\ &\qquad
+\sum_{r\in\{1,3,5\}}\biggl(\biggl\lf\frac{j+a_r-a_0-1}p\biggr\rf-\biggl\lf\frac{j+b_r-a_0-1}p\biggr\rf-\biggl\lf\frac{a_r-b_r}p\biggr\rf\biggr).
\label{eq130}
\end{align}
Combining the two cases together we conclude that
\begin{align}
\ord_pB(\ba,\bb), \; \ord_p d_{m_1}^3d_{m_2}C(\ba,\bb)
&\ge\min\{\Omega_1(\ba,\bb;p),\Omega_2(\ba,\bb;p)\}
\nonumber
\end{align}
for primes $p>\sqrt{b_0-a_0-2}$, where the quantities $\Omega_1$ and $\Omega_2$ are defined in \eqref{eq128} and~\eqref{eq130}.

When we choose
\begin{equation}
\begin{gathered}
a_j=\alpha_jn+1 \quad\text{for}\; j=0,1,\dots,6,
\\
b_0=\beta_0n+3 \quad\text{and}\quad
b_j=\beta_jn+1 \quad\text{for}\; j=1,\dots,6,
\end{gathered}
\label{eq132}
\end{equation}
for some positive set of integer directions $(\balpha,\bbeta)$,
then computing $\Omega_1,\Omega_2$ reduces to the computation of the minima $\omega_1^*(x)$ and $\omega_2^*(x)$ of functions
\begin{align}
\omega_1(x,y,z)
&=\lf(\beta_0-\alpha_0)x\rf-\lf y+z-\alpha_0x\rf-\lf\beta_0x-(y+z)\rf
\nonumber\\ &\qquad
+\sum_{r\in\{2,4,6\}}\bigl(\lf y-\beta_rx\rf-\lf y-\alpha_rx\rf-\lf(\alpha_r-\beta_r)x\rf\bigr)
\nonumber\\ &\qquad
+\sum_{r\in\{1,3,5\}}\bigl(\lf z-\beta_rx\rf-\lf z-\alpha_rx\rf-\lf(\alpha_r-\beta_r)x\rf\bigr)
\nonumber
\\ \intertext{and}
\omega_2(x,y)
&=\sum_{r\in\{2,4,6\}}\bigl(\lf y-\beta_rx\rf-\lf y-\alpha_rx\rf-\lf(\alpha_r-\beta_r)x\rf\bigr)
\nonumber\\ &\qquad
+\sum_{r\in\{1,3,5\}}\bigl(\lf y+(\alpha_r-\alpha_0)x\rf-\lf y+(\beta_r-\alpha_0)x\rf-\lf(\alpha_r-\beta_r)x\rf\bigr)
\nonumber
\end{align}
over $y,z$ and over $y$, respectively. Indeed,
$$
\Omega_1(\ba,\bb;p)=\omega_1\biggl(\frac np,\frac{j-1}p,\frac{l-1}p\biggr)
\quad\text{and}\quad
\Omega_2(\ba,\bb;p)=\omega_2\biggl(\frac np,\frac{j-1}p\biggr)
$$
in the settings above. This means that
\begin{equation}
\ord_pB(\ba,\bb), \; 4+\ord_pC(\ba,\bb)
\ge\min\biggl\{\omega_1^*\biggl(\frac np\biggr),\omega_2^*\biggl(\frac np\biggr)\biggr\}
\label{eq135}
\end{equation}
for primes $p>\sqrt{(\beta_0-\alpha_0)n}$.

Notice that the functions $\omega_1$ and $\omega_2$ (hence their minima) are $1$-periodic in each variable,
so it is sufficient to compute them on the intervals $[0,1)$.
In the most symmetric case
$$
\alpha_0=\beta_1=\dots=\beta_6=0, \quad
\alpha_1=\dots=\alpha_6=1 \quad\text{and}\quad \beta_0=3
$$
we already get (by droping the four non-negative terms in both $\omega_1$ and $\omega_2$)
\begin{align*}
&
\lf3x\rf-\lf y+z\rf-\lf3x-(y+z)\rf+(\lf y\rf-\lf y-x\rf-\lf x\rf)+(\lf z\rf-\lf z-x\rf-\lf x\rf)
\\&\quad
=\lf3x\rf-\lf y+z\rf-\lf3x-(y+z)\rf-\lf y-x\rf-\lf z-x\rf
\ge1 \quad\text{for}\; x\in[\tfrac23,1)
\end{align*}
and
\begin{align*}
&
(\lf y\rf-\lf y-x\rf-\lf x\rf)+(\lf y+x\rf-\lf y\rf-\lf x\rf)
\\ &\quad
=\lf y+x\rf-\lf y-x\rf
\ge1 \quad\text{for}\; x\in[\tfrac12,1).
\end{align*}
Indeed, when $x\in[\tfrac23,1)$, the first inequality follows from
$$
\lf3x\rf-\lf y+z\rf-\lf3x-(y+z)\rf\ge0 \quad\text{and}\quad \lf y-x\rf+\lf z-x\rf\le-1
$$
if either $y<x$ or $z<x$; otherwise, $\frac23\le x\le y<1$ and $\frac23\le x\le z<1$ imply
\begin{gather*}
\lf3x\rf=2, \quad \lf y+z\rf=1, \quad \lf3x-(y+z)\rf=\lf x-(y-x)-(z-x)\rf=0
\\
\text{and}\quad \lf y-x\rf=\lf z-x\rf=0.
\end{gather*}
A proof of the second inequality, when $x\in[\frac12,1)$, makes use of $\lf y-x\rf=-1$, $\lf y+x\rf\ge0$ if $0\le y<\frac12$, and of $\lf y+x\rf=1$, $\lf y-x\rf\le0$ if $\frac12\le y<1$.
The two inequalities together mean that the quantity $\Phi_n$ from Section~\ref{sec:sym} is divisible by~\eqref{eq9}.

It looks quite plausible (though we do not possess any proof of this) that we always have $\omega_2^*(x)\ge\omega_1^*(x)$ except for possibly finitely many rational points on the interval $[0,1)$. (Notice that $\omega_1(x,y,\alpha_0x-y)$ coincides with $\omega_2(x,y)$ apart from finitely many rational lines crossing the square $[0,1)^2$.)

\medskip
Now assume that the linear forms $G(\ba,\bb)$ originate from the forms $F(\bh)$ of Section~\ref{sec:old} and condition \eqref{eq30} written as
\begin{equation}
2\beta_0+\alpha_0=\alpha_1+\alpha_2+\alpha_3+\alpha_4+\alpha_5+\alpha_6
\nonumber
\end{equation}
holds. In this case we can scale all the parameters in \eqref{eq28} to discuss the set $\be n$ instead, where
$$
\begin{gathered}
e_{0j}=\alpha_j, \; \ol e_{0j}=\alpha_j-\alpha_0 \quad\text{for}\; 1\le j\le6,
\\
e_{jk}=\beta_0-\alpha_j-\alpha_k \quad\text{for}\; 1\le j<k\le6,
\end{gathered}
$$
and record the related quantities by $H(\be n)=B(\be n)\zeta(4)-C(\be n)\in\mathbb Z\zeta(4)+\mathbb Q$, where $n=0,1,2,\dots$\,.
The discussion above (see \eqref{eq135}) implies that
$$
\ord_pB(\be n), \, 4+\ord_pC(\be n)
\ge\omega^*\biggl(\be;\frac np\biggr)
\quad\text{for primes}\; p>\sqrt{(\beta_0-\alpha_0)n},
$$
where
$$
\omega^*(\be;x)=\min\{\omega_1^*(\be;x),\omega_2^*(\be;x)\},
$$
hence also
$$
\ord_pB(\fg\be n), \, 4+\ord_pC(\fg\be n)
\ge\omega^*\biggl(\fg\be;\frac np\biggr)
\quad\text{for primes}\; p>\sqrt{(\beta_0-\alpha_0)n} \;\text{and}\; \fg\in\fG,
$$
where $\fg\be$ denotes the image of the multiset $\be$ under the action of $\fg\in\fG$.
At the same time,
$$
\frac{B(\be n)}{\Pi(\be n)}=\frac{B(\fg\be n)}{\Pi(\fg\be n)}
\quad\text{and}\quad
\frac{C(\be n)}{\Pi(\be n)}=\frac{C(\fg\be n)}{\Pi(\fg\be n)}
$$
for all $\fg\in\fG$, in view of the invariance of $H(\be n)/\Pi(\be n)$ under the action of $\fG$
(and of the irrationality of $\zeta(4)$). This implies that
\begin{align*}
\ord_pB(\be n), \, 4+\ord_pC(\be n)
&\ge\ord_p\biggl(\frac{\Pi(\be n)}{\Pi(\fg\be n)}\,\omega^*\biggl(\fg\be;\frac np\biggr)\biggr)
\\
&=\sum_{e\in\cE}\biggl(\biggl\lf\frac{en}p\biggr\rf-\biggl\lf\frac{\fg en}p\biggr\rf\biggr)
+\omega^*\biggl(\fg\be;\frac np\biggr)
\end{align*}
for primes $p>\sqrt{(\beta_0-\alpha_0)n}$ and all $\fg\in\fG$,
hence
\begin{equation}
\ord_pB(\be n), \, 4+\ord_pC(\be n)
\ge\omega\biggl(\be;\frac np\biggr)
\quad\text{for primes}\; p>\sqrt{(\beta_0-\alpha_0)n},
\nonumber
\end{equation}
where
\begin{equation}
\omega(\be;x)
=\max_{\fg\in\fG}\biggl(\sum_{e\in\cE}(\lf ex\rf-\lf\fg ex\rf)+\omega^*(\fg\be;x)\biggr).
\label{eq138}
\end{equation}
The maximum can be restricted to distinct representatives modulo the group of $(\ba,\bb)$-trivial permutations.

\section({One concrete example of irrationality measure for \003\266(4)})%
{One concrete example of irrationality measure for $\zeta(4)$}
\label{sec:ex}

In the notation of Section~\ref{sec:old} we take
$$
h_0=\eta_0n+2, \; h_{-1}=\eta_{-1}n+2, \; h_1=\eta_1n+1, \; \dots, \; h_6=\eta_6n+1
$$
with
$$
\Beta=(\eta_0,\eta_{-1};\eta_1,\dots,\eta_6)=(68,57;22,23,24,25,26,27).
$$
If we set $F_n=F(\bh)=G(\ba,\bb)=u_n\zeta(4)-v_n$ then the asymptotics of $F_n$ and $u_n$ as $n\to\infty$ is computed
with the help of \cite[Proposition~1]{Zu03} (adapted here to address a a slightly different normalisation of $F(\bh)$):
$$
C_0=-\lim_{n\to\infty}\frac{\log|F_n|}n=36.47011287\dots
\quad\text{and}\quad
C_1=\lim_{n\to\infty}\frac{\log|u_n|}n=106.34774225\dots\,.
$$

The above choice of $\bh$ translates the form $F_n=F(\bh)$ from \eqref{eq10} into $G(\ba,\bb)$ from \eqref{eq11} with the parameters \eqref{eq132} as follows:
\begin{equation}
\balpha=(11; 22, 23, 24, 25, 26, 27), \quad \bbeta=(68; 0, 0, 4, 7, 11, 11).
\label{bab}
\end{equation}
The denominator of $v_n=C(\ba,\bb)$ in \eqref{eq26} is $d_{21n}^3d_{23n}^{\vphantom1}$.
The following table lists 31 out of 120 representatives under the action of group $\fG$ on \eqref{bab} modulo the trivial $(\ba,\bb)$-action,
only those that contribute to the computation of the corresponding function $\omega(x)=\omega(\be;x)$ in~\eqref{eq138}:

\begin{center}
{\scriptsize
\begin{tabular}{|r|c|}
\hline
1 & (68; 22, 23, 24, 25, 26, 27) \\
2 & (68; 22, 23, 24, 25, 27, 26) \\
3 & (68; 22, 23, 24, 26, 25, 27) \\
4 & (68; 22, 23, 25, 24, 26, 27) \\
5 & (68; 22, 23, 25, 24, 27, 26) \\
6 & (68; 22, 23, 26, 24, 27, 25) \\
7 & (68; 22, 24, 23, 25, 26, 27) \\
8 & (68; 22, 24, 23, 25, 27, 26) \\
9 & (68; 22, 25, 23, 26, 24, 27) \\
10 & (67; 21, 22, 23, 25, 26, 27) \\
& \\
\hline
\end{tabular}
\hfil
\begin{tabular}{|r|c|}
\hline
11 & (67; 21, 22, 23, 25, 27, 26) \\
12 & (67; 21, 22, 23, 26, 25, 27) \\
13 & (67; 21, 22, 25, 23, 26, 27) \\
14 & (66; 20, 21, 23, 24, 26, 27) \\
15 & (66; 20, 21, 23, 24, 27, 26) \\
16 & (66; 20, 21, 23, 26, 24, 27) \\
17 & (66; 20, 21, 24, 23, 27, 26) \\
18 & (65; 19, 20, 23, 24, 25, 27) \\
19 & (65; 19, 20, 23, 24, 27, 25) \\
20 & (65; 19, 20, 23, 25, 24, 27) \\
21 & (65; 19, 20, 24, 23, 27, 25) \\
\hline
\end{tabular}
\hfil
\begin{tabular}{|r|c|}
\hline
22 & (65; 19, 21, 22, 23, 26, 27) \\
23 & (65; 19, 21, 22, 23, 27, 26) \\
24 & (65; 19, 21, 23, 22, 26, 27) \\
25 & (65; 19, 21, 23, 22, 27, 26) \\
26 & (65; 19, 21, 26, 22, 27, 23) \\
27 & (65; 19, 22, 21, 23, 26, 27) \\
28 & (65; 19, 23, 20, 24, 27, 25) \\
29 & (64; 18, 19, 23, 25, 24, 26) \\
30 & (64; 19, 20, 21, 22, 27, 26) \\
31 & (64; 19, 21, 20, 22, 26, 27) \\
& \\
\hline
\end{tabular}
\hfil
}
\end{center}

\noindent
Here we give the representatives in the format $(\beta_0;\alpha_1,\dots,\alpha_6)=(\eta_0;\eta_1,\dots,\eta_6)$;
all other parameters are completely determined by the data.

Then
\begin{align*}
\omega(x)=0 \quad\text{if}\;
x&\in\bigl[0,\tfrac2{57}\bigr)
\cup\bigl[\tfrac1{15},\tfrac2{27}\bigr)
\cup\bigl[\tfrac7{16},\tfrac{25}{57}\bigr)
\cup\bigl[\tfrac{11}{23},\tfrac{28}{57}\bigr)
\cup\bigl[\tfrac8{15},\tfrac7{13}\bigr)
\cup\bigl[\tfrac{22}{23},\tfrac{56}{57}\bigr),
\displaybreak[2]\\
\omega(x)=1 \quad\text{if}\;
x&\in\bigl[\tfrac2{57},\tfrac1{27}\bigr)\perm{1}
\cup\bigl[\tfrac1{22},\tfrac1{21}\bigr)\perm{1}
\cup\bigl[\tfrac1{17},\tfrac1{16}\bigr)\perm{1}
\cup\bigl[\tfrac1{16},\tfrac1{15}\bigr)\perm{3}
\cup\bigl[\tfrac2{27},\tfrac1{13}\bigr)\perm{1}
\\ &\quad
\cup\bigl[\tfrac2{15},\tfrac3{22}\bigr)\perm{1}
\cup\bigl[\tfrac3{22},\tfrac8{57}\bigr)\perm{4}
\cup\bigl[\tfrac17,\tfrac4{27}\bigr)\perm{1}
\cup\bigl[\tfrac5{23},\tfrac29\bigr)\perm{1}
\cup\bigl[\tfrac5{21},\tfrac{14}{57}\bigr)\perm{1}
\\ &\quad
\cup\bigl[\tfrac4{15},\tfrac3{11}\bigr)\perm{1}
\cup\bigl[\tfrac27,\tfrac7{24}\bigr)\perm{1}
\cup\bigl[\tfrac5{17},\tfrac8{27}\bigr)\perm{1}
\cup\bigl[\tfrac7{20},\tfrac{20}{57}\bigr)\perm{2}
\cup\bigl[\tfrac6{17},\tfrac5{14}\bigr)\perm{1}
\\ &\quad
\cup\bigl[\tfrac25,\tfrac{23}{57}\bigr)\perm{2}
\cup\bigl[\tfrac{10}{23},\tfrac7{16}\bigr)\perm{1}
\cup\bigl[\tfrac{25}{57},\tfrac{11}{25}\bigr)\perm{1}
\cup\bigl[\tfrac{10}{21},\tfrac{11}{23}\bigr)\perm{1}
\cup\bigl[\tfrac{28}{57},\tfrac12\bigr)\perm{1}
\\ &\quad
\cup\bigl[\tfrac9{17},\tfrac8{15}\bigr)\perm{1}
\cup\bigl[\tfrac7{13},\tfrac{31}{57}\bigr)\perm{1}
\cup\bigl[\tfrac{14}{23},\tfrac{11}{18}\bigr)\perm{1}
\cup\bigl[\tfrac{11}{18},\tfrac{35}{57}\bigr)\perm{4}
\cup\bigl[\tfrac{12}{17},\tfrac{41}{57}\bigr)\perm{1}
\\ &\quad
\cup\bigl[\tfrac{17}{23},\tfrac34\bigr)\perm{1}
\cup\bigl[\tfrac{13}{17},\tfrac{10}{13}\bigr)\perm{1}
\cup\bigl[\tfrac45,\tfrac{46}{57}\bigr)\perm{1}
\cup\bigl[\tfrac{19}{23},\tfrac56\bigr)\perm{1}
\cup\bigl[\tfrac{20}{23},\tfrac78\bigr)\perm{1}
\\ &\quad
\cup\bigl[\tfrac78,\tfrac{50}{57}\bigr)\perm{2}
\cup\bigl[\tfrac{19}{21},\tfrac{10}{11}\bigr)\perm{1}
\cup\bigl[\tfrac{10}{11},\tfrac{52}{57}\bigr)\perm{7}
\cup\bigl[\tfrac{21}{23},\tfrac{11}{12}\bigr)\perm{1}
\cup\bigl[\tfrac{20}{21},\tfrac{21}{22}\bigr)\perm{1}
\\ &\quad
\cup\bigl[\tfrac{21}{22},\tfrac{22}{23}\bigr)\perm{7}
\cup\bigl[\tfrac{56}{57},1\bigr)\perm{1},
\displaybreak[2]\\
\omega(x)=2 \quad\text{if}\;
x&\in\bigl[\tfrac1{27},\tfrac1{26}\bigr)\perm{1}
\cup\bigl[\tfrac1{23},\tfrac1{22}\bigr)\perm{1}
\cup\bigl[\tfrac1{19},\tfrac1{18}\bigr)\perm{1}
\cup\bigl[\tfrac1{18},\tfrac1{17}\bigr)\perm{2}
\cup\bigl[\tfrac1{13},\tfrac2{25}\bigr)\perm{1}
\\ &\quad
\cup\bigl[\tfrac2{19},\tfrac19\bigr)\perm{1}
\cup\bigl[\tfrac19,\tfrac3{26}\bigr)\perm{2}
\cup\bigl[\tfrac2{17},\tfrac7{57}\bigr)\perm{1}
\cup\bigl[\tfrac3{23},\tfrac2{15}\bigr)\perm{3}
\cup\bigl[\tfrac8{57},\tfrac17\bigr)\perm{4}
\displaybreak[1]\\ &\quad
\cup\bigl[\tfrac4{27},\tfrac3{20}\bigr)\perm{20}
\cup\bigl[\tfrac3{20},\tfrac2{13}\bigr)\perm{3}
\cup\bigl[\tfrac3{19},\tfrac16\bigr)\perm{1}
\cup\bigl[\tfrac4{23},\tfrac{10}{57}\bigr)\perm{11}
\cup\bigl[\tfrac3{17},\tfrac2{11}\bigr)\perm{1}
\displaybreak[1]\\ &\quad
\cup\bigl[\tfrac2{11},\tfrac5{27}\bigr)\perm{5}
\cup\bigl[\tfrac4{21},\tfrac{11}{57}\bigr)\perm{6}
\cup\bigl[\tfrac{11}{57},\tfrac15\bigr)\perm{1}
\cup\bigl[\tfrac15,\tfrac5{24}\bigr)\perm{3}
\cup\bigl[\tfrac4{19},\tfrac5{23}\bigr)\perm{1}
\displaybreak[1]\\ &\quad
\cup\bigl[\tfrac29,\tfrac5{22}\bigr)\perm{2}
\cup\bigl[\tfrac5{22},\tfrac{13}{57}\bigr)\perm{4}
\cup\bigl[\tfrac{13}{57},\tfrac3{13}\bigr)\perm{2}
\cup\bigl[\tfrac4{17},\tfrac5{21}\bigr)\perm{4}
\cup\bigl[\tfrac{14}{57},\tfrac14\bigr)\perm{1}
\displaybreak[1]\\ &\quad
\cup\bigl[\tfrac5{19},\tfrac4{15}\bigr)\perm{3}
\cup\bigl[\tfrac3{11},\tfrac5{18}\bigr)\perm{14}
\cup\bigl[\tfrac5{18},\tfrac{16}{57}\bigr)\perm{4}
\cup\bigl[\tfrac{16}{57},\tfrac27\bigr)\perm{2}
\cup\bigl[\tfrac7{24},\tfrac5{17}\bigr)\perm{2}
\displaybreak[1]\\ &\quad
\cup\bigl[\tfrac8{27},\tfrac{17}{57}\bigr)\perm{1}
\cup\bigl[\tfrac7{23},\tfrac4{13}\bigr)\perm{3}
\cup\bigl[\tfrac6{19},\tfrac7{22}\bigr)\perm{1}
\cup\bigl[\tfrac7{22},\tfrac8{25}\bigr)\perm{4}
\cup\bigl[\tfrac8{25},\tfrac9{26}\bigr)\perm{1}
\displaybreak[1]\\ &\quad
\cup\bigl[\tfrac8{23},\tfrac7{20}\bigr)\perm{21}
\cup\bigl[\tfrac{20}{57},\tfrac6{17}\bigr)\perm{2}
\cup\bigl[\tfrac5{14},\tfrac4{11}\bigr)\perm{3}
\cup\bigl[\tfrac7{19},\tfrac38\bigr)\perm{1}
\cup\bigl[\tfrac8{21},\tfrac5{13}\bigr)\perm{1}
\displaybreak[1]\\ &\quad
\cup\bigl[\tfrac9{23},\tfrac25\bigr)\perm{1}
\cup\bigl[\tfrac{23}{57},\tfrac{11}{27}\bigr)\perm{2}
\cup\bigl[\tfrac7{17},\tfrac5{12}\bigr)\perm{3}
\cup\bigl[\tfrac8{19},\tfrac{10}{23}\bigr)\perm{1}
\cup\bigl[\tfrac{11}{25},\tfrac49\bigr)\perm{10}
\displaybreak[1]\\ &\quad
\cup\bigl[\tfrac49,\tfrac{26}{57}\bigr)\perm{2}
\cup\bigl[\tfrac9{19},\tfrac{10}{21}\bigr)\perm{4}
\cup\bigl[\tfrac12,\tfrac{29}{57}\bigr)\perm{2}
\cup\bigl[\tfrac{10}{19},\tfrac9{17}\bigr)\perm{2}
\cup\bigl[\tfrac{31}{57},\tfrac6{11}\bigr)\perm{1}
\displaybreak[1]\\ &\quad
\cup\bigl[\tfrac6{11},\tfrac{11}{20}\bigr)\perm{10}
\cup\bigl[\tfrac{11}{20},\tfrac59\bigr)\perm{5}
\cup\bigl[\tfrac47,\tfrac{15}{26}\bigr)\perm{2}
\cup\bigl[\tfrac{11}{19},\tfrac7{12}\bigr)\perm{1}
\cup\bigl[\tfrac{10}{17},\tfrac{13}{22}\bigr)\perm{1}
\displaybreak[1]\\ &\quad
\cup\bigl[\tfrac{13}{22},\tfrac{16}{27}\bigr)\perm{6}
\cup\bigl[\tfrac{16}{27},\tfrac{34}{57}\bigr)\perm{1}
\cup\bigl[\tfrac35,\tfrac{14}{23}\bigr)\perm{1}
\cup\bigl[\tfrac{35}{57},\tfrac8{13}\bigr)\perm{4}
\cup\bigl[\tfrac{13}{21},\tfrac{17}{27}\bigr)\perm{2}
\displaybreak[1]\\ &\quad
\cup\bigl[\tfrac{12}{19},\tfrac7{11}\bigr)\perm{1}
\cup\bigl[\tfrac7{11},\tfrac{16}{25}\bigr)\perm{4}
\cup\bigl[\tfrac{11}{17},\tfrac{37}{57}\bigr)\perm{1}
\cup\bigl[\tfrac{15}{23},\tfrac{17}{26}\bigr)\perm{3}
\cup\bigl[\tfrac{17}{26},\tfrac{17}{25}\bigr)\perm{1}
\displaybreak[1]\\ &\quad
\cup\bigl[\tfrac{13}{19},\tfrac{11}{16}\bigr)\perm{1}
\cup\bigl[\tfrac{11}{16},\tfrac9{13}\bigr)\perm{4}
\cup\bigl[\tfrac{16}{23},\tfrac7{10}\bigr)\perm{14}
\cup\bigl[\tfrac7{10},\tfrac{40}{57}\bigr)\perm{4}
\cup\bigl[\tfrac{40}{57},\tfrac{12}{17}\bigr)\perm{2}
\displaybreak[1]\\ &\quad
\cup\bigl[\tfrac{41}{57},\tfrac{18}{25}\bigr)\perm{1}
\cup\bigl[\tfrac{14}{19},\tfrac{17}{23}\bigr)\perm{7}
\cup\bigl[\tfrac34,\tfrac{43}{57}\bigr)\perm{4}
\cup\bigl[\tfrac{43}{57},\tfrac{19}{25}\bigr)\perm{2}
\cup\bigl[\tfrac{16}{21},\tfrac{13}{17}\bigr)\perm{2}
\displaybreak[1]\\ &\quad
\cup\bigl[\tfrac{10}{13},\tfrac{44}{57}\bigr)\perm{1}
\cup\bigl[\tfrac{18}{23},\tfrac{11}{14}\bigr)\perm{3}
\cup\bigl[\tfrac{15}{19},\tfrac45\bigr)\perm{1}
\cup\bigl[\tfrac{46}{57},\tfrac{21}{26}\bigr)\perm{1}
\cup\bigl[\tfrac{17}{21},\tfrac{13}{16}\bigr)\perm{1}
\displaybreak[1]\\ &\quad
\cup\bigl[\tfrac{13}{16},\tfrac9{11}\bigr)\perm{2}
\cup\bigl[\tfrac{14}{17},\tfrac{47}{57}\bigr)\perm{8}
\cup\bigl[\tfrac{47}{57},\tfrac{19}{23}\bigr)\perm{1}
\cup\bigl[\tfrac56,\tfrac{21}{25}\bigr)\perm{3}
\cup\bigl[\tfrac{16}{19},\tfrac{11}{13}\bigr)\perm{1}
\displaybreak[1]\\ &\quad
\cup\bigl[\tfrac{13}{15},\tfrac{20}{23}\bigr)\perm{7}
\cup\bigl[\tfrac{50}{57},\tfrac{22}{25}\bigr)\perm{2}
\cup\bigl[\tfrac{15}{17},\tfrac89\bigr)\perm{16}
\cup\bigl[\tfrac{17}{19},\tfrac9{10}\bigr)\perm{1}
\cup\bigl[\tfrac9{10},\tfrac{19}{21}\bigr)\perm{3}
\\ &\quad
\cup\bigl[\tfrac{52}{57},\tfrac{21}{23}\bigr)\perm{7}
\cup\bigl[\tfrac{11}{12},\tfrac{23}{25}\bigr)\perm{27}
\cup\bigl[\tfrac{23}{25},\tfrac{12}{13}\bigr)\perm{1}
\cup\bigl[\tfrac{18}{19},\tfrac{19}{20}\bigr)\perm{1}
\cup\bigl[\tfrac{19}{20},\tfrac{20}{21}\bigr)\perm{3},
\displaybreak[2]\\
\omega(x)=3 \quad\text{if}\;
x&\in\bigl[\tfrac1{26},\tfrac1{25}\bigr)\perm{1}
\cup\bigl[\tfrac1{21},\tfrac1{20}\bigr)\perm{19}
\cup\bigl[\tfrac1{20},\tfrac1{19}\bigr)\perm{23}
\cup\bigl[\tfrac2{25},\tfrac1{12}\bigr)\perm{14}
\cup\bigl[\tfrac2{23},\tfrac5{57}\bigr)\perm{26}
\\ &\quad
\cup\bigl[\tfrac5{57},\tfrac1{11}\bigr)\perm{13}
\cup\bigl[\tfrac2{21},\tfrac1{10}\bigr)\perm{19}
\cup\bigl[\tfrac1{10},\tfrac2{19}\bigr)\perm{23}
\cup\bigl[\tfrac3{26},\tfrac2{17}\bigr)\perm{2}
\cup\bigl[\tfrac7{57},\tfrac18\bigr)\perm{1}
\displaybreak[1]\\ &\quad
\cup\bigl[\tfrac18,\tfrac3{23}\bigr)\perm{3}
\cup\bigl[\tfrac2{13},\tfrac3{19}\bigr)\perm{19}
\cup\bigl[\tfrac16,\tfrac4{23}\bigr)\perm{7}
\cup\bigl[\tfrac{10}{57},\tfrac3{17}\bigr)\perm{11}
\cup\bigl[\tfrac5{27},\tfrac3{16}\bigr)\perm{30}
\displaybreak[1]\\ &\quad
\cup\bigl[\tfrac3{16},\tfrac4{21}\bigr)\perm{14}
\cup\bigl[\tfrac5{24},\tfrac4{19}\bigr)\perm{23}
\cup\bigl[\tfrac3{13},\tfrac4{17}\bigr)\perm{4}
\cup\bigl[\tfrac14,\tfrac7{27}\bigr)\perm{2}
\cup\bigl[\tfrac6{23},\tfrac5{19}\bigr)\perm{23}
\displaybreak[1]\\ &\quad
\cup\bigl[\tfrac{17}{57},\tfrac3{10}\bigr)\perm{1}
\cup\bigl[\tfrac3{10},\tfrac7{23}\bigr)\perm{3}
\cup\bigl[\tfrac4{13},\tfrac5{16}\bigr)\perm{19}
\cup\bigl[\tfrac5{16},\tfrac6{19}\bigr)\perm{21}
\cup\bigl[\tfrac9{26},\tfrac8{23}\bigr)\perm{19}
\displaybreak[1]\\ &\quad
\cup\bigl[\tfrac4{11},\tfrac7{19}\bigr)\perm{19}
\cup\bigl[\tfrac38,\tfrac8{21}\bigr)\perm{10}
\cup\bigl[\tfrac5{13},\tfrac{22}{57}\bigr)\perm{14}
\cup\bigl[\tfrac{22}{57},\tfrac7{18}\bigr)\perm{1}
\cup\bigl[\tfrac7{18},\tfrac9{23}\bigr)\perm{2}
\displaybreak[1]\\ &\quad
\cup\bigl[\tfrac{11}{27},\tfrac9{22}\bigr)\perm{2}
\cup\bigl[\tfrac9{22},\tfrac7{17}\bigr)\perm{15}
\cup\bigl[\tfrac5{12},\tfrac8{19}\bigr)\perm{23}
\cup\bigl[\tfrac{26}{57},\tfrac{11}{24}\bigr)\perm{2}
\cup\bigl[\tfrac8{17},\tfrac9{19}\bigr)\perm{19}
\displaybreak[1]\\ &\quad
\cup\bigl[\tfrac{29}{57},\tfrac{14}{27}\bigr)\perm{2}
\cup\bigl[\tfrac{11}{21},\tfrac{10}{19}\bigr)\perm{21}
\cup\bigl[\tfrac59,\tfrac{32}{57}\bigr)\perm{5}
\cup\bigl[\tfrac{13}{23},\tfrac47\bigr)\perm{17}
\cup\bigl[\tfrac{15}{26},\tfrac{11}{19}\bigr)\perm{19}
\displaybreak[1]\\ &\quad
\cup\bigl[\tfrac7{12},\tfrac{10}{17}\bigr)\perm{10}
\cup\bigl[\tfrac{34}{57},\tfrac35\bigr)\perm{1}
\cup\bigl[\tfrac8{13},\tfrac{13}{21}\bigr)\perm{15}
\cup\bigl[\tfrac{17}{27},\tfrac{12}{19}\bigr)\perm{21}
\cup\bigl[\tfrac{16}{25},\tfrac9{14}\bigr)\perm{14}
\displaybreak[1]\\ &\quad
\cup\bigl[\tfrac9{14},\tfrac{11}{17}\bigr)\perm{10}
\cup\bigl[\tfrac{37}{57},\tfrac{13}{20}\bigr)\perm{1}
\cup\bigl[\tfrac{13}{20},\tfrac{15}{23}\bigr)\perm{8}
\cup\bigl[\tfrac{17}{25},\tfrac{15}{22}\bigr)\perm{24}
\cup\bigl[\tfrac{15}{22},\tfrac{13}{19}\bigr)\perm{19}
\displaybreak[1]\\ &\quad
\cup\bigl[\tfrac9{13},\tfrac{16}{23}\bigr)\perm{14}
\cup\bigl[\tfrac{18}{25},\tfrac{13}{18}\bigr)\perm{10}
\cup\bigl[\tfrac{13}{18},\tfrac8{11}\bigr)\perm{2}
\cup\bigl[\tfrac8{11},\tfrac{19}{26}\bigr)\perm{6}
\cup\bigl[\tfrac{11}{15},\tfrac{14}{19}\bigr)\perm{19}
\displaybreak[1]\\ &\quad
\cup\bigl[\tfrac{19}{25},\tfrac{16}{21}\bigr)\perm{4}
\cup\bigl[\tfrac{44}{57},\tfrac{17}{22}\bigr)\perm{1}
\cup\bigl[\tfrac{17}{22},\tfrac79\bigr)\perm{8}
\cup\bigl[\tfrac79,\tfrac{18}{23}\bigr)\perm{3}
\cup\bigl[\tfrac{11}{14},\tfrac{15}{19}\bigr)\perm{19}
\displaybreak[1]\\ &\quad
\cup\bigl[\tfrac{21}{26},\tfrac{17}{21}\bigr)\perm{12}
\cup\bigl[\tfrac9{11},\tfrac{14}{17}\bigr)\perm{14}
\cup\bigl[\tfrac{21}{25},\tfrac{16}{19}\bigr)\perm{23}
\cup\bigl[\tfrac{11}{13},\tfrac{17}{20}\bigr)\perm{14}
\cup\bigl[\tfrac{17}{20},\tfrac67\bigr)\perm{3}
\displaybreak[1]\\ &\quad
\cup\bigl[\tfrac67,\tfrac{49}{57}\bigr)\perm{9}
\cup\bigl[\tfrac{49}{57},\tfrac{19}{22}\bigr)\perm{1}
\cup\bigl[\tfrac{19}{22},\tfrac{13}{15}\bigr)\perm{7}
\cup\bigl[\tfrac{22}{25},\tfrac{15}{17}\bigr)\perm{14}
\cup\bigl[\tfrac89,\tfrac{17}{19}\bigr)\perm{19}
\\ &\quad
\cup\bigl[\tfrac{12}{13},\tfrac{13}{14}\bigr)\perm{14}
\cup\bigl[\tfrac{16}{17},\tfrac{17}{18}\bigr)\perm{23}
\cup\bigl[\tfrac{17}{18},\tfrac{18}{19}\bigr)\perm{19},
\displaybreak[2]\\
\omega(x)=4 \quad\text{if}\;
x&\in\bigl[\tfrac1{25},\tfrac1{24}\bigr)\perm{14}
\cup\bigl[\tfrac1{24},\tfrac1{23}\bigr)\perm{10}
\cup\bigl[\tfrac1{12},\tfrac2{23}\bigr)\perm{24}
\cup\bigl[\tfrac1{11},\tfrac2{21}\bigr)\perm{31}
\cup\bigl[\tfrac7{27},\tfrac6{23}\bigr)\perm{25}
\\ &\quad
\cup\bigl[\tfrac{11}{24},\tfrac6{13}\bigr)\perm{11}
\cup\bigl[\tfrac7{15},\tfrac8{17}\bigr)\perm{18}
\cup\bigl[\tfrac{14}{27},\tfrac{13}{25}\bigr)\perm{11}
\cup\bigl[\tfrac{12}{23},\tfrac{11}{21}\bigr)\perm{22}
\cup\bigl[\tfrac{32}{57},\tfrac9{16}\bigr)\perm{5}
\\ &\quad
\cup\bigl[\tfrac9{16},\tfrac{13}{23}\bigr)\perm{12}
\cup\bigl[\tfrac{19}{26},\tfrac{11}{15}\bigr)\perm{19}
\cup\bigl[\tfrac{13}{14},\tfrac{53}{57}\bigr)\perm{19}
\cup\bigl[\tfrac{14}{15},\tfrac{15}{16}\bigr)\perm{28}
\cup\bigl[\tfrac{15}{16},\tfrac{16}{17}\bigr)\perm{29},
\displaybreak[2]\\
\omega(x)=5 \quad\text{if}\;
x&\in\bigl[\tfrac6{13},\tfrac7{15}\bigr)\perm{18}
\cup\bigl[\tfrac{13}{25},\tfrac{12}{23}\bigr)\perm{22}
\cup\bigl[\tfrac{53}{57},\tfrac{14}{15}\bigr)\perm{19},
\end{align*}
where the notation $[a,b)\perm{N}$ means that the maximum in~\eqref{eq138} is attained on the $N$-th representative.
Denoting
$$
\Phi_n=\prod_{p>\sqrt{57n}}p^{\omega(n/p)}
$$
we conclude that $\Phi_n^{-1}u_n\in\mathbb Z$ and $\Phi_n^{-1}d_{21n}^3d_{23n}v_n\in\mathbb Z$;
in other words,
$$
\Phi_n^{-1}d_{21n}^3d_{23n}F_n\in\mathbb Z\zeta(4)+\mathbb Z
\quad\text{for}\; n=1,2,\dots\,.
$$
At the same time the asymptotics of $\Phi_n^{-1}d_{21n}^3d_{23n}$ is controlled by the prime number theorem:
$$
C_2=\lim_{n\to\infty}\frac{\log(\Phi_n^{-1}d_{21n}^3d_{23n})}n
=3\cdot21+23-\int_0^1\omega(x)\,\d\psi(x)
=25.05460171\dots,
$$
where $\psi(x)$ is the logarithmic derivative of the gamma function. Now \cite[Proposition~3]{Zu03} applies to imply
that the irrationality exponent of $\zeta(4)$ is bounded above by
$$
\frac{C_0+C_1}{C_0-C_2}=12.51085940\dots\,.
$$

\medskip
Finally, we point out that the general family of rational approximations to $\zeta(4)$ from Section~\ref{sec:gen}
is only exploited here when it is linked to the old approximations reviewed in Section~\ref{sec:old}.
A reason behind this is mainly an easy access to the asymptotic behaviour of the corresponding forms $G(\ba,\bb)$ and their coefficients $B(\ba,\bb)$.
One may hope to get a better control of general approximations from Section~\ref{sec:gen} by covering analytic aspects of the 12-parametric family there,
however this will not necessarily lead to (significantly) better arithmetic consequences.


\end{document}